\documentclass[10pt]{article}
\usepackage{graphicx}
\usepackage{amsmath}
\usepackage{amsfonts}
\usepackage{amssymb}
\newtheorem{theorem}{Theorem}
\newtheorem{corollary}[theorem]{Corollary}
\newtheorem{lemma}[theorem]{Lemma}
\newtheorem{proposition}[theorem]{Proposition}
\newtheorem{definition}[theorem]{Definition}
\newtheorem{remark}{Remark}
\newtheorem{example}{Example}
\newcommand{\N}{\mathbb{N}}
\newcommand{\R}{\mathbb{R}}
\newcommand{\dom}{\mathop\mathrm{\rm dom}}
\newcommand{\dist}{\mbox{\rm dist\,}}
\newcommand{\crit}{\mathop\mathrm{\rm crit}}
\newcommand{\graph}{\mbox{Graph\,}}

\newcommand{\diag}{\mbox{\rm diag\,}}
\newcommand{\rang}{\mbox{\rm rank\,}}
\newcommand{\s}{\mathbb{S}}

\newcommand{\Frechet}{\hat{\partial}}

\newcommand{\argmin}{\mbox{\rm argmin}\;}
\newcommand{\xk}{x^k}

\def\eps{\varepsilon}

\begin{document}

\begin{center}
{\LARGE 
  Proximal alternating minimization and projection methods for nonconvex problems.  
  An approach based on the Kurdyka-{\L}ojasiewicz inequality } 
\footnote{The first three authors acknowledge the support of the French ANR 
under grant ANR-08-BLAN-0294-03.}

\vspace{1cm}
\bf{Hedy ATTOUCH
  \footnote{Institut de Math\'ematiques et de Mod\'elisation de Montpellier, 
  UMR CNRS 5149, CC 51, Universit\'e de Montpellier II, Place Eug\`ene 
  Bataillon, 34095 Montpellier cedex 5, France (attouch@math.univ-montp2.fr,
  redont@math.univ-montp2.fr).},}  
\bf{ J\'er\^ome BOLTE
  \footnote{UPMC Paris 06,  Equipe Combinatoire et Optimisation and Inria Saclay (CMAP, Polytechnique) , UMR 7090, case 189, 
  Universit\'e Pierre et Marie Curie, 4 Place Jussieu, 75252 Paris cedex 05
  (bolte@math.jussieu.fr).},}
\bf{ Patrick REDONT $^2$,} 
\bf{ Antoine SOUBEYRAN
  \footnote{GREQAM, UMR CNRS 6579, Universit\'e de la M\'editerran\'ee, 13290 
  Les Milles, France (antoine.soubeyran@univmed.fr).}.}

\bigskip
\end{center}

\noindent\textbf{Abstract.} 
\small
We study the convergence properties of an alternating proximal minimization 
algorithm for nonconvex structured functions of the type:
$L(x,y)=f(x)+Q(x,y)+g(y)$, where $f:\R^n\rightarrow\R\cup\{+\infty\}$ and 
$g:\R^m\rightarrow\R\cup\{+\infty\}$ are proper lower semicontinuous functions, 
and $Q:\R^n\times\R^m\rightarrow \R$ is a smooth $C^1$ function which couples 
the variables $x$ and $y$. The algorithm can be viewed as a {\em proximal} 
regularization of the usual Gauss-Seidel method to minimize $L$.

We work in a nonconvex setting, just assuming that the function $L$ satisfies 
the Kurdyka-\L ojasiewicz inequality. An entire section illustrates the 
relevancy of such an assumption by giving examples ranging from semialgebraic 
geometry to ``metrically regular" problems. 

Our main result can be stated as follows: If L has the 
Kurdyka-\L ojasiewicz property, then  each bounded sequence generated by the algorithm 
 converges to a critical point of $L$. 
This result is completed by the study of the convergence rate of the 
algorithm, which depends on the geometrical properties of the function $L$ 
around its critical points. When specialized to $Q(x,y)=\left\|x-y\right\|^2$
and to $f$, $g$ indicator functions, the algorithm is an alternating projection 
mehod (a variant of Von Neumann's) that converges for a wide class of sets 
including semialgebraic and tame sets, transverse smooth manifolds or sets 
with ``regular ''intersection. In order to illustrate our results with concrete problems, we provide a 
convergent proximal reweighted $\ell^1$ algorithm for compressive sensing and an application to 
rank reduction problems.
 

\medskip


\normalsize
\section{Introduction}

\noindent
{\bf Presentation of the algorithm.} In this paper, we will be concerned with the convergence analysis of 
alternating minimization algorithms for (nonconvex) functions 
$L:\R^n\times\R^m\rightarrow\R\cup\{+\infty\}$ of the following type:

$(\mathcal H)\quad\left\{\rule{0em}{3em}\right.$
\parbox{39em}{
  $L(x,y)=f(x)+Q(x,y)+g(y)$,\\
  $f:\R^n\rightarrow\R\cup\{+\infty\}$, $g:\R^m\rightarrow\R\cup\{+\infty\}$ 
    are proper lower semicontinuous,\\
  $Q:\R^n\times\R^m\rightarrow\R$ is a $C^1$ function,\\
  $\nabla Q$ is Lipschitz continuous on bounded subsets of $\R^n\times\R^m$. 
}

Assumption $(\mathcal H)$ will be needed throughout the paper.

We aim at finding critical points of 

\begin{equation}
\label{Pmin} L(x,y)=f(x)+Q(x,y)+g(y)
\end{equation}

\noindent and possibly solve the corresponding minimization problem. 

The specific structure of $L$ allows in particular to tackle problems of the 
form
\begin{equation}
\label{Pmin2}\min\{f(z)+g(z):z\in\R^n\}.
\end{equation} 
It suffices indeed to set $L_\rho(x,y)=f(x)+\frac{\rho}{2}||x-y||^2+g(y)$, 
$\rho$ being a positive penalization (or relaxation) parameter, and to 
minimize $L_\rho$ over $\R^n\times\R^n$. Feasibility problems involving two 
closed sets are particular cases of the above problem: just specialize $f$ and 
$g$ to be the indicator functions. (In a convex setting, the square of the 
euclidean distance may be replaced by a Bregman distance, see \cite{BCN}.) 

Minimizing  the sum of simple functions or finding a common point to a 
collection of closed sets is a very active field of research with applications 
in approximation theory \cite{Johannis}, image reconstruction 
\cite{CombWajs,donoho},  statistics \cite{CT,GrubP}, 
partial differential equations and optimal control \cite{PL.Lions,WW}. A good 
reference for problems involving convex instances is  \cite{CombWajs}: many 
examples coming from signal processing problems are shown to be rewritable 
as~(\ref{Pmin}).  

The specificity of our approach is twofold. First, we work in a nonconvex 
setting, just assuming that the function $L$ satisfies the 
Kurdyka-\L ojasiewicz inequality, see \cite{Loja63,Loja93,Kurdyka98}. As it 
has been established recently in \cite{BDL04,BDL04b,BDLS}, this assumption is 
satisfied by a wide class of nonsmooth functions called functions definable in 
an o-minimal structure (see Section \ref{S:omin}). Semialgebraic functions 
and (globally) subanalytic functions are for instance definable in their 
respective classes. 
 
Secondly, we rely on a new class of alternating minimization algorithms with 
costs to move which has recently been introduced in \cite{ARS} (see also \cite{AS}) and which has 
proved to be a flexible tool \cite{abrs} permitting to handle general coupling 
functions $Q(x,y)$ (for example $Q(x,y)=\|Ax-By\|^2$, with $A$, $B$ linear 
operators): 

$(x_0,y_0)\in\R^n\times\R^m$ given, 
$(x_k,y_k)\rightarrow(x_{k+1},y_k)\rightarrow(x_{k+1},y_{k+1})$
\begin{equation}\label{algo0}
\left\{
\begin{aligned}
 x_{k+1} & \in\mbox{argmin}\left\{L(u,y_k)+\frac{1}{2\lambda_k}\|u-x_k\|^2:u\in\R^n\right\},\\
 y_{k+1} & \in \mbox{argmin}\left\{L(x_{k+1},v)+\frac{1}{2\mu_k}\|v-y_k\|^2:v\in \R^m\right\}.
\end{aligned}\right.
\end{equation}

The above algorithm can be viewed as a \textit{proximal} regularization of a 
two block Gauss-Seidel method for minimizing $L$:
$$\left\{
\begin{array}{l}
x_{k+1}\in\mbox{argmin}\{L(u,y_k):u\in\R^n\}\\
\rule{0pt}{16pt}
y_{k+1}\in\mbox{argmin}\{L(x_{k+1},v):v\in\R^m\}.
\end{array}\right.$$
Some general results for Gauss-Seidel method, also known as coordinate descent 
method, can be found for instance in \cite{Aus,Ber}; block coordinate methods 
for nonsmooth and nonconvex functions have been investigated by many authors 
(see \cite{Tseng} and references therein). However very few general results 
ensure that the sequence $(x_k,y_k)$ converges to a global minimizer, even 
for strictly convex functions. An important fact concerning our approach is 
that the convergence of algorithm (\ref{algo0}) works for any stepsizes 
$\lambda_k,\mu_k$ greater than a fixed positive parameter which can be chosen 
arbitrarily large. For such parameters the algorithm is very close to a 
coordinate descent method. On the other hand, when the stepsizes are not too 
large, the method is an alternating gradient-like method.

\medskip

\noindent
{\bf Kurdyka-\L ojasiewicz inequalities and tame geometry.} Before describing 
and illustrating our convergence results, let us recall some important facts 
that have motivated our mathematical approach. 

In his pioneering work on real-analytic functions \cite{Loja63,Loja93},  
\L ojasiewicz provided the basic ingredient, the so-called ``\L ojasiewicz 
inequality'', that allows to derive the convergence of the bounded 
trajectories of the steepest descent equation to critical points. Given a 
real-analytic function $f:\R^n\rightarrow \R$ and a critical point $a\in\R^n,$ 
the \L ojasiewicz inequality asserts that there exists some 
$\theta\in\left[\frac{1}{2},1\right)$ such that the function 
$|f-f(a)|^\theta\left\|\nabla f\right\|^{-1}$ remains bounded around $a$. 
Similar results have been developed for discrete gradient methods (see 
\cite{AMA2004}) and nonsmooth subanalytic functions (see \cite{BDL04,BDL04b,attbol}). 
  
In the last decades powerful advances relying on an axiomatized approach of 
real-semialgebraic/real-analytic geometry have allowed to set up a general 
theory in which the basic objects enjoy the same qualitative properties as 
semialgebraic sets and functions \cite{Dries,Shiota,Wilkie}. In such a 
framework the central concept is that of o-minimal structure over $\R$. 
Basic results are recalled and illustrated in Section \ref{S:omin}. Following 
van den Dries \cite{Dries}, functions and sets belonging to such structures 
are called definable or {\em tame}~(\footnote{The word tame actually 
corresponds to a slight generalization of definable objects.}). Extensions of 
the \L ojasiewicz inequality to definable functions and applications to their 
gradient vector fields have been obtained by Kurdyka \cite{Kurdyka98}, while 
nonsmooth versions have been developed in \cite{BDLS}. The corresponding
generalized inequality is here called the Kurdyka-\L ojasiewicz inequality  
(see Definition \ref{D:Kurdyk}, Section \ref{S:Kurdyk}).

An important motivation for designing alternating algorithms for tame problems 
relies on these generalized \L ojasiewicz inequalities but also on the 
following facts: 

-- tame sets and functions allow to model many problems both smooth and 
nonsmooth: most of standard subsets of matrices are semialgebraic (symmetric 
positive semidefinite matrices, orthogonal groups, Stiefel manifolds, constant 
rank matrices) and most criteria involve piecewise polynomial and analytic 
functions.

-- tameness is a highly stable ``concept'': finite unions or finite 
intersections of tame sets are tame, compositions of tame mappings are tame, 
subdifferentials of tame functions are tame (see the paragraph after Definition 
\ref{D:omin} and also \cite{BDL04,BDLS}).

\medskip

\noindent
{\bf Convergence results, tractability of the algorithm and applications. }
One of our central result (Theorem \ref{T.bounded}) can be stated as follows: 
Assume that $L$ has the Kurdyka-\L ojasiewicz property at each point. Then
either the sequence $(x_k,y_k)$ converges to infinity, or the trajectory has a 
finite length and, as a consequence, converges to a critical point of $L$. 
This result is completed by the study of the convergence rate of the 
algorithm, which depends on the geometrical properties of the function $L$ 
around its critical points (namely the \L ojasiewicz exponent). 

%The relevancy of an algorithm in Optimization relies on several crucial 
%features among which verifiability of the assumptions, computability, 
%efficiency are probably the most important. However, it is also of primary 
%importance, to identify simple mechanisms that  rule convergence aspects. 
%When one or some of these key properties have been enlighted,  it becomes
%easier to target other specific properties that might impact the usual 
%triplet verifiability/computability/efficiency and raise in turn new 
%theoretical issues.
%
\smallskip
Let us now give some insight into the interest of algorithm~(\ref{algo0}). 

\smallskip
{\em Verifying} that a nonsmooth function has the Kurdyka-\L ojasiewicz 
property {\em at each point} is often a very easy task. For instance check the 
semialgebraicity or the analyticity of functions $f$, $g$, $Q$ and apply the 
nonsmooth \L ojasiewicz inequalities provided in \cite{BDL04,BDLS}. Section 
\ref{S.Appli} below shows the preeminence of this inequality and its links to 
prominent concepts in Optimization: 

-- uniformly convex functions and convex functions enjoying  growth conditions;

-- metric regularity and constraint qualification;

-- semialgebraic and definable functions. 

Specific examples related to feasibility problems are provided, they involve 
(possibly tangent) real-analytic manifolds, transverse manifolds (see 
\cite{mallew}), semialgebraic sets or more generally tame sets.

Our results also apply to {\em tame convex problems} which are easily 
identifiable in practice and which provide in turn an important field of 
applications (see for instance \cite{CombWajs}). Many convergence results of 
various types are available under convexity assumptions \cite{BB,CombWajs}. 
Our results seem new; when $L$ is a tame closed convex function and has at 
least a minimizer, then the sequences generated by algorithm (\ref{algo0}) 
converge. Besides, when $L$ is semialgebraic or globally subanalytic, 
convergence rates are necessarily of the form  $O(\frac{1}{k^s})$ with $s>0$.
 
{\em Computability} of the proximal iterates deserves some explanation. For 
simplicity take $Q(x,y)=\frac{1}{2}||x-y||^2$. A first general remark is the 
following: when, for instance, $f$ is locally convex up to a square, an 
adequate choice of the stepsize $\lambda_k$ makes the computation of the first 
step of algorithm \ref{algo0} a {\em convex}, hence tractable,  problem. 
Computational issues concerning these aspects of the implementation of 
 approximate proximal points are given by Hare-Sagastiz\'abal in \cite{Hare}; besides, 
the Nesterov optimal gradient algorithm provides a simple tool for efficiently 
solving convex problems (\cite{Nesterov}). Some nonconvex cases are easily 
computable and amount to {\em explicit} or standard computations: projections 
onto spheres or onto ellipsoids, constant rank matrices, or even one real 
variable second order equations. A good reference for some of these aspects is 
\cite{mallew}. In order to provide a more realistic and flexible tool for solving real-world problems,  it would be  natural to consider inexact versions of algorithm (\ref{algo0}) (see \cite{IPS,CP04} for some work in that direction); this subject is out of the scope of the present paper but it is a matter for future research.

\medskip
Several applications involving nonconvex aspects are provided (see Section 
\ref{S:omin}): rank reduction of correlation matrices, compressive sensing 
with  nonconvex ``norms''. For this last case our algorithm provides a 
regularized version of the reweighted $l^1$ minimization algorithm of 
Cand\`es-Wakin-Boyd (see Example~\ref{compressed}). To the best of our 
knowledge no convergence results are known for the reweighted algorithms. The
regularized version we provide converges even when small constant stepsizes 
are chosen in the implementation of (\ref{algo0}).   

In Section \ref{vonN}, when specializing  algorithm  (\ref{algo0}) to 
indicator functions, we obtain an  alternating projection algorithm (recently 
introduced in \cite{ARS}), which can be seen as a proximal regularization of 
the von Neumann algorithm (see \cite{Johannis}). Being given two closed 
subsets $C,D$ of $\R^n$ the algorithm reads
\begin{equation*}
\left\{\begin{aligned}
x_{k+1} & \in  P_C\left(\frac{\lambda_k^{-1}x_k+y_k}{\lambda_k^{-1}+1}\right)\\
y_{k+1} & \in  P_D\left(\frac{\mu_k^{-1}y_k+x_{k+1}}{\mu_k^{-1}+1}\right),
\end{aligned}\right.
\end{equation*}
where $P_C, P_D :\R^n\rightrightarrows \R^n$ are the projection  mappings onto 
$C$ and $D$. The convergence of the sequences $(x_k),(y_k)$ is obtained for a 
wide class of sets ranging from semialgebraic or definable sets to transverse 
manifolds, or more generally to sets with a regular intersection. A part of 
this result  is inspired by the recent work of Lewis and Malick on transverse 
manifolds \cite{mallew} (and also \cite{LLM}), in which similar results were derived.

\bigskip

The paper is organized as follows: Section 2 is devoted to recalling some 
elementary facts of nonsmooth analysis. This allows us to obtain in Section 3 
some first elementary properties of the alternating proximal minimization 
algorithm, and then to establish our main theoretical results. A last section 
is devoted to examples and applications: various classes of functions 
satisfying the Kurdyka-Lojasiewicz property are provided and specific examples 
for which computations can be effectively performed are given. 

\section{Elementary facts of nonsmooth analysis}
  
The Euclidean scalar product of $\R^n$ and its corresponding  norm are 
respectively denoted by $\langle\cdot,\cdot\rangle$ and $||\cdot||$.
Some general references for nonsmooth analysis are \cite{Rock98,Morduk}.

If $F:\R^n\rightrightarrows\R^m$ is a point-to-set mapping its {\em graph} is 
defined by
$$\graph F:=\{(x,y)\in\R^n\times\R^m:y\in F(x)\}.$$
Similarly the graph of a real-extended-valued function 
$f:\R^n\rightarrow \R\cup\{+\infty\}$ is defined by
$$\graph f:=\{(x,s)\in\R^n\times\R: s=f(x)\}.$$

Let us recall a few definitions concerning subdifferential calculus. 

\begin{definition}{\rm (\cite{Rock98})   Let $f:\R^n\rightarrow \R\cup\{+\infty\}$ be a proper lower semicontinuous function. 

(i) The domain of $f$ is defined and denoted by $\dom f:=\{x\in \R^n: f(x)<+\infty\}$.

(ii) For each $x\in \dom f$, the Fr\'echet subdifferential of $f$ at $x$, written $\Frechet f(x)$, is
 the set of vectors $x^*\in \R^n$ which satisfy 

$$\liminf_{\begin{array}{l}y\neq x \\ y\rightarrow x\end{array}} \frac{1}{\|x-y\|}[f(y)-f(x)-\langle x^*,y-x\rangle]\geq0.$$

If $x\notin\dom f$, then $\Frechet f(x)=\emptyset$.

(iii) The limiting-subdifferential (\cite{Morduk}), or simply the 
subdifferential for short, of $f$ at $x\in\dom f$, written $\partial f(x)$, is 
defined as follows 
$$\partial f(x):=
  \{x^*\in\R^n:
    \exists x_n\rightarrow x,\;f(x_n)\rightarrow f(x),\;
      x_n^*\in\Frechet f(x_n)\rightarrow x^*
  \}.
$$
}
\end{definition} 

\begin{remark}\label{subdiff} 
{\rm (a) The above definition implies that $\Frechet f(x)\subset\partial f(x)$ 
for each $x\in\R^n$, where the first set is convex and closed while the second 
one is closed \cite[th. 8.6 p. 302]{Rock98}.\\
(b)(Closedness of $\partial f$) Let $(x_k,x_k^*)\in\graph\partial f$ be a 
sequence that converges to $(x,x^*)$. By the very definition of 
$\partial f(x)$, if $f(x_k)$ converges to $f(x)$ then 
$(x,x^*)\in\graph\partial f$.\\
(c) A necessary (but not sufficient) condition for $x\in\R^n$ to be a 
minimizer of $f$ is
\begin{equation}\label{crit}
\partial f(x)\ni 0.
\end{equation}

A point that satisfies (\ref{crit}) is called {\em limiting-critical} or simply critical. 
The set 
of critical points of $f$ is denoted by $\crit f$.}
\end{remark}

If $K$ is a subset of $\R^n$ and $x$ is any point in $\R^n$, we set
$$\dist(x,K)=\inf\{\|x-z\|:z\in K\}.$$
Recall that if $K$ is empty we have $\dist(x,K)=+\infty$ for all $x\in\R^n$. 
Note also that for any real-extended-valued function $f$ on $\R^n$ and any 
$x\in\R^n$, $\dist(0,\partial f(x))=\inf\{\|x^*\|:x^*\in\partial f(x)\}$.

\begin{lemma}\label{semicontinuity} Let $f:\R^n\rightarrow \R\cup\{+\infty\}$ be a proper
 lower semicontinuous function. 
Let $\bar x\in\dom f$ be a noncritical point of $f$. Then there exists $c>0$ 
such that
$$\|x-\bar x\|+\|f(x)-f(\bar x)\|<c  \ \Longrightarrow  \  \dist(0,\partial f(x))\geq c.$$
\end{lemma}
{\bf Proof}. On the contrary, there would exist a sequence $(c_k)$ with $c_k >0$,
$c_k\rightarrow 0$, and a sequence $(x_k)$ with 
$\|x_k-\bar x\|+\|f(x_k)-f(\bar x)\|<c_k$ and $\dist(0,\partial f(x_k))<c_k$. 
The latter inequality implies the existence of some $x_k^*\in\partial f(x_k)$ 
with $\|x_k^*\|<c_k$. Owing to the closedness of $\partial f$ we would then 
have $0\in\partial f(\bar x)$, a contradiction.$\hfill\Box$

\medskip
\noindent
{\bf Partial subdifferentiation} Let 
$L:\R^n\times\R^m\rightarrow\R\cup\{+\infty\}$ be a lower semicontinuous 
function. When fixing $y$   in $\R^m$, the subdifferential of the 
function $L(\cdot,y)$ at $u$ is denoted by $\partial_x L(u,y)$. Similarly, when fixing $x$
in $\R^n$, one 
can define the partial subdifferentiation with respect to the variable $y$. The 
corresponding operator is denoted by $\partial_y L(x,\cdot).$ 

The following result, though elementary, is central to the paper.
\begin{proposition}\label{sous-diff}
Let $L$ satisfy $(\mathcal H)$. Then for all 
$(x,y)\in\dom L=\dom f\times\dom g$ we have 
$$
\partial L(x,y)
=\left\{\partial f(x)+\nabla_x Q(x,y)\right\}\times
 \left\{\partial g(y)+\nabla_y Q(x,y)\right\}
=\partial_x L(x,y)\times\partial_y L(x,y).
$$
\end{proposition}
{\bf Proof.} 
Observe first that we have
$\partial L(x,y)=\partial(f(x)+g(y))+\nabla Q(x,y)$, since $Q$ is continuously 
differentiable (\cite[8.8(c) Exercice, p. 304]{Rock98}). Further, the 
subdifferential calculus for separable functions yields 
(\cite[10.5 Proposition, p. 426]{Rock98}) 
$\partial(f(x)+g(y))=\partial f(x)\times\partial g(y)$. Hence the first 
equality.

Invoking once more (\cite[8.8(c) Exercice]{Rock98}) yields the second 
equality.$\hfill\Box$

\medskip
\noindent
{\bf Normal cones, indicator functions and projections} 

If $C$ is a closed subset of $\R^n$ we denote by $\delta_C$ its indicator function, {\it i.e.} for all $x\in \R^n$
 we set  
$$\delta_C(x)=\left\{
\begin{array}{ll}
0 & \mbox{ if }x\in C,\\
+\infty & \mbox{otherwise.}
\end{array}\right.
$$ 

The projection on $C$, written $P_C$, is the following {\em point-to-set} 
mapping:
$$P_C: \left\{
\begin{array}{lll}
\R^n & \rightrightarrows & \R^n\\
x & \rightarrow & P_C(x):=\argmin\{\|x-z\|:z\in C\}.
\end{array}\right.$$
When $C$ is nonempty, the closedness of $C$ and the compactness of the closed 
unit ball of $\R^n$ imply that $P_C(x)$ is nonempty for all $x$ in $\R^n$.

\begin{definition}{\rm (Normal cone) 
Let $C$ be a nonempty closed subset of $\R^n$.\\ 
(i) For any $x\in C$ the Fr\'echet normal cone to $C$ at $x$ is defined by
$$\hat{N}_C(x)=\{v\in\R^n:\langle v,y-x\rangle\leq o(x-y), \ y\in C\}.$$
When $x\notin C$ we set $N_C(x)=\emptyset$.\\
(ii) The (limiting) normal cone to $C$ at $x\in C$ is denoted by $N_C(x)$ and 
is defined by
$$v\in N_C(x)\Leftrightarrow
  \exists x_k\in C, x_k\rightarrow x,\,
  \exists v_k\in\hat{N}_C(x_k), v_k\rightarrow v.$$}
\end{definition}

\begin{remark} {\rm (a) For $x\in C$ the cone $N_C(x)$ is closed but not necessarily convex.\\
(b) An elementary but important fact about normal cone and subdifferential is the following 
$$\partial \delta_C=N_C.$$

}\end{remark}

\section{Alternating proximal minimization algorithms}

\subsection{Convergence to a critical value}\label{algo}

Let $L$ satisfy $(\mathcal H)$. Being given $(x_0,y_0)\in\R^n\times\R^m$,
recall that the alternating discrete dynamical system we are to study is of 
the form: 
$(x_k,y_k)\rightarrow(x_{k+1},y_k)\rightarrow(x_{k+1},y_{k+1})$
\begin{flushright}
\makebox[0pt]
  {\rule{140pt}{0pt}\raisebox{-5pt}{$\left\{\rule{0pt}{23pt}\right.$}}
\begin{minipage}{400pt}
\begin{eqnarray}
x_{k+1}&\in&\argmin \{L(u,y_k)+\frac{1}{2\lambda_k}
\|u-x_k\|^2: u\in \R^n\}\label{prox1}\\
y_{k+1}&\in&\argmin \{L(x_{k+1},v)+\frac{1}{2\mu_k}\|v-y_k\|^2: v\in \R^m\},\label{prox2}
\end{eqnarray}
\end{minipage}
\end{flushright}
where $(\lambda_k)_{k\in\N},(\mu_k)_{k\in\N}$ are positive sequences. 

\noindent
We make the following standing assumption concerning (\ref{prox1}), 
(\ref{prox2}): 
$$(\mathcal H_1) \left\{
\begin{array}{l}\quad\inf_{\R^n\times\R^m}L>-\infty,\\
\quad\mbox{the function }L(\cdot,y_0)\mbox{ is proper,}\\
\quad\mbox{for some positive }r_-<r_+\mbox{ the sequences of stepsizes } 
\lambda_k,\;\mu_k\mbox{ belong to }(r_-,r_+)\mbox{ for all }k\geq 0.
\end{array}\right.$$

The next lemma, especially point (iii), is of constant use in the sequel.
\begin{lemma}\label{basic}
Under assumptions $(\mathcal H)$, $(\mathcal H_1)$, the sequences $(x_k)$, 
$(y_k)$ are correctly defined. Moreover\\
(i) The following estimate holds
\begin{equation}\label{descent}
  L(x_k,y_k)
    +\frac{1}{2\lambda_{k-1}}\|x_k-x_{k-1}\|^2
    +\frac{1}{2\mu_{k-1}}\|y_k-y_{k-1}\|^2
  \leq L(x_{k-1},y_{k-1})\quad\forall k\geq 1;
\end{equation}
hence $L(x_k,y_k)$ does not increase.\\
(ii) 
$$\sum_{k=1}^\infty\left(\|x_k-x_{k-1}\|^2+\|y_k-y_{k-1}\|^2\right)<+\infty;$$
hence $\lim(\|x_k-x_{k-1}\|+\|y_k-y_{k-1}\|)=0$.\\
(iii) For $k\geq 1$  define 
$(x_k^*,y_k^*)=
  (\nabla_x Q(x_k,y_k) -\nabla_x Q(x_k,y_{k-1}),0)
    -\left(\frac{1}{\lambda_{k-1}}(x_k -x_{k-1}),
      \frac{1}{\mu_{k-1}}(y_k -y_{k-1})\right)$;
we have 
\begin{equation}\label{zetoile}
(x_k^*,y_k^*)\in\partial L(x_k,y_k).
\end{equation}
For all bounded subsequence $(x_{k'},y_{k'})$ of $(x_k,y_k)$ we have 
$(x_{k'}^*,y_{k'}^*)\rightarrow 0,\;k'\rightarrow+\infty$, hence\linebreak 
$\dist(0,\partial L(x_{k'},y_{k'}))\rightarrow 0,\;k'\rightarrow+\infty$.
\end{lemma}
\noindent
{\bf Proof.} Since $\inf L>-\infty$, $(\cal H)$ implies that for any $r>0$, 
$(\bar u,\bar v)\in\R^n\times\R^m$ the functions 
$u\rightarrow L(u,\bar v)+\frac{1}{2r}\|u-\bar u\|^2$ and 
$v\rightarrow L(\bar u,v)+\frac{1}{2r}\|v-\bar v\|^2$ are  coercive. An 
elementary induction ensures then that the sequences are well defined and 
that (i) and (ii) hold for all integer $k\geq 1$. 

By the very definition of $x_k$, and Remark \ref{subdiff} c), $0$ must lie in 
the subdifferential at point $x_k$ of the function 
$\xi\mapsto\frac{1}{2\lambda_{k-1}}\|\xi-{k-1}\|^2+L(\xi,y_{k-1})$ which is 
equal to $\frac{1}{\lambda_{k-1}}(\xi-x_{k-1})+\partial_x L(x_k,y_{k-1})$ 
since the function $\xi\mapsto\frac{1}{2\lambda_{k-1}}\|\xi-x_{k-1}\|^2$ is 
smooth. Hence 
\begin{equation}\label{prox'1}
0\in\frac{1}{\lambda_{k-1}}(x_k -x_{k-1})+\partial_x L(x_k,y_{k-1}),  
 \ \forall k\geq 1.
\end{equation}
And similarly
\begin{equation}\label{prox'2}
0\in\frac{1}{\mu_{k-1}}(y_k -y_{k-1})+\partial_y L(x_k,y_k),  
 \ \forall k\geq 1.
\end{equation}
Due to the structure of $L$ we have  
$\partial_x L(x_k,y_{k-1})=\partial f(x_k)+\nabla_x Q(x_k,y_{k-1})$ and 
$\partial_y L(x_k,y_k)=\partial g(y_k)+\nabla_y Q(x_k,y_k)$. Hence we may 
write with (\ref{prox'1}, \ref{prox'2}) 
$$
\begin{array}{c}
-\frac{1}{\lambda_{k-1}}(x_k -x_{k-1})
  -(\nabla_x Q(x_k,y_{k-1})-\nabla_x Q(x_k,y_k),0)
\in\partial f(x_k)+\nabla_x Q(x_k,y_k);\\
\\
-\frac{1}{\mu_{k-1}}(y_k -y_{k-1})\in\partial g(y_k)+\nabla_y Q(x_k,y_k).
\end{array}
$$
This yields (\ref{zetoile}) with proposition \ref{sous-diff}.
 
If $(x_{k'},y_{k'})$ is a bounded sequence, then so is $(x_{k'},y_{k'-1})$, 
and, by (ii), $(x_{k'},y_{k'})-(x_{k'},y_{k'-1})$ vanishes as 
$k'\rightarrow\infty$. The uniform continuity of $\nabla_x Q$ on bounded 
subsets then yields the last point of (iii).$\hfill\Box$
 
\begin{remark}{\rm (a) Without additional assumptions, like for instance the 
convexity of $L$, the sequence $(x_k,y_k)$ is not {\it a priori} uniquely 
defined.\\
(b) Note also that the result remains valid if $L$ is bounded from below by an affine function.}
\end{remark}

The next proposition, especially points (ii)(iii), gives the first convergence 
results about sequences generated by (\ref{prox1},\ref{prox2}). Theorems 
\ref{T.bounded} and \ref{T.rate} below make the convergence properties much 
more precise.
\begin{proposition}\label{GEN} Assume that $(\mathcal H)$, $(\mathcal H_1)$ 
hold. Let $(x_k,y_k)$ be a sequence complying with (\ref{prox1}) and (\ref{prox2}). 
Let $\omega(x_0,y_0)$ denote the set (possibly empty) of its limit points. 
Then\\
(i) If $(x_k,y_k)$ is bounded, then $\omega (x_0,y_0)$ is a nonempty compact 
connected set and
$$d((x_k,y_k),\omega(x_0,y_0))\rightarrow 0\mbox{ as }k\rightarrow+\infty,$$
(ii) $\omega(x_0,y_0)\subset\crit L$,\\
(iii) $L$ is finite and constant on $\omega(x_0,y_0)$, equal to 
$\inf_{k\in\N}L(x_k,y_k)=\lim_{k\rightarrow+\infty}L(x_k,y_k)$.
\end{proposition}
{\bf Proof.} 
Item (i) follows by using $\|x_k-x_{k-1}\|+\|y_k-y_{k-1}\|\rightarrow 0$ 
together with some classical properties of sequences in $\R^n$.\\
(ii) By the very definition of $(x_k,y_k)$ ($k\geq 1$) we have
$$
\begin{array}{c}
L(x_k,y_{k-1})+\frac{1}{2\lambda_{k-1}}\|x_k-x_{k-1}\|^2
  \leq L(\xi,y_{k-1})+\frac{1}{2\lambda_{k-1}}\|\xi-x_{k-1}\|^2,\quad
  \forall\xi\in\R^m\\
L(x_k,y_k)+\frac{1}{2\mu_{k-1}}\|y_k-y_{k-1}\|^2
  \leq L(x_k,\eta)+\frac{1}{2\mu_{k-1}}\|\eta-y_{k-1}\|^2,\quad
  \forall\eta\in\R^n.\\
\end{array}
$$
Due to the special form of $L$ and to $0<r_-\leq\lambda_{k-1}\leq r_+$ and 
$0<r_-\leq\mu_{k-1}\leq r_+$, we have
\begin{eqnarray}
&  f(x_k)+Q(x_k,y_{k-1})+\frac{1}{2r_+}\|x_k-x_{k-1}\|^2
  \leq f(\xi)+Q(\xi,y_{k-1})+\frac{1}{2r_-}\|\xi-x_{k-1}\|^2,\quad
  \forall\xi\in\R^m, & \label{toto}\\
& g(y_k)+Q(x_k,y_k)+\frac{1}{2r_+}\|y_k-y_{k-1}\|^2
  \leq g(\eta)+Q(x_k,\eta)+\frac{1}{2r_-}\|\eta-y_{k-1}\|^2,\quad
  \forall\eta\in\R^n. & \label{titi}
\end{eqnarray}
Let $(\bar x,\bar y)$ be a point in $\omega(x_0,y_0)$; there exists a 
subsequence $(x_{k'},y_{k'})$ of $(x_k,y_k)$ converging to $(\bar x,\bar y)$. 
Since $\|x_k-x_{k-1}\|+\|y_k-y_{k-1}\|\rightarrow 0$ we deduce from 
(\ref{toto})
$$
\liminf f(x_{k'})+Q(\bar x,\bar y)\leq f(\xi)+Q(\xi,\bar y) +\frac{1}{2r_-}\|\xi-\bar x \|^2      \quad 
  \forall\xi\in\R^m.$$
In particular for $\xi=\bar x$ we obtain
$$\liminf f(x_{k'})\leq f(\bar x).$$
Since $f$ is lower semicontinuous we get further
$$\liminf f(x_{k'})=f(\bar x).$$
There is no loss of generality in assuming that the whole sequence $f(x_{k'})$ 
converges to $f(\bar x)$
\begin{equation}\label{convfk}
\lim f(x_{k'})=f(\bar x).
\end{equation}
Similarly, using (\ref{titi}), we may assume that 
$\lim g(y_{k'})=g(\bar y).$ 
Now, as $Q$ is continuous, we have \linebreak  
$Q(x_{k'},y_{k'})\rightarrow Q(\bar x,\bar y)$ and hence  
$L(x_{k'},y_{k'})\rightarrow L(\bar x,\bar y)$. But with lemma 
\ref{basic}(iii) and using the same notation , we have $(x_{k'}^*,y_{k'}^*)\in\partial L(x_{k'},y_{k'})$ and 
$(x_{k'}^*,y_{k'}^*)\rightarrow 0$. Owing to the closedness properties of $\partial L$, we finally 
obtain $0\in\partial L(\bar x,\bar y)$.

(iii) For any point $(\bar x,\bar y)\in\omega(x_0,y_0)$ we have just seen that 
there is a subsequence $(x_{k'},y_{k'})$ with 
$L(x_{k'},y_{k'})\rightarrow L(\bar x,\bar y)$. Since the sequence 
$L(x_k,y_k)$ is nonincreasing, we have $ L(\bar x,\bar y)=\inf L(x_k,y_k)$ 
independent of $(\bar x,\bar y)$.$\hfill\Box$

\subsection{Kurdyka-\L ojasiewicz inequality}\label{S:Kurdyk}

Let $f:\R^n\rightarrow \R\cup\{+\infty\}$ be a proper lower semicontinuous
function. For $-\infty<\eta_1<\eta_2\leq+\infty$, let us set 
$$[\eta_1<f<\eta_2]=\left\{x\in \R^n: \eta_1<f(x)<\eta_2\right\}.$$
\begin{definition}[Kurdyka-\L ojasiewicz property]\label{D:Kurdyk}
The function $f$ is said to have the Kurdyka-\L ojasie\-wicz property at 
$\bar x\in\dom\partial f$ if there exist $\eta\in (0,+\infty]$, a 
neighborhood (\footnote{See Remark (\ref{r1}) (c)}) $U$ of $\bar x$ and a 
continuous {\em concave} function $\varphi:[0,\eta)\rightarrow\R_+$ such that:

{\bf -} $\varphi(0)=0$,

{\bf -} $\varphi$ is $C^1$ on $(0,\eta)$,

{\bf -} for all $s\in (0,\eta)$, $\varphi'(s)>0$,

{\bf -} and for all $x$ in $U\cap[f(\bar x)<f<f(\bar x)+\eta]$, the 
Kurdyka-\L ojasie\-wicz inequality holds
\begin{equation}\label{Loja}
\varphi'(f(x)-f(\bar x))\,\dist(0,\partial f(x))\geq 1.
\end{equation}
\end{definition}

\begin{remark}\label{r1}
{\rm (a) S. \L ojasiewicz proved in 1963 \cite{Loja63} that real-analytic 
functions satisfy an inequality of the above type with $\varphi(s)=s^{1-\theta}$ 
where $\theta\in[\frac{1}{2},1)$. A nice proof of this result can be found in 
the monograph \cite{Denk}. In a recent paper, Kurdyka \cite{Kurdyka98} has 
extended this result to differentiable functions definable in an o-minimal 
structure (see  Section~\ref{S:omin}). More recently Bolte, Daniilidis, Lewis and 
Shiota have extended the Kurdyka-Lojasiewicz inequality to nonsmooth functions 
in the subanalytic and o-minimal setting \cite{BDL04,BDL04b,BDLS}; see also   
\cite{BDLM}.\\
The concavity assumption imposed on the function $\varphi$ does not explicitly 
belong to the usual formulation of the Kurdyka-Lojasiewicz inequality. However
 the 
examples of the next section show that the inequality holds in 
many instances with a 
concave $\varphi$.\\ 
(b) A proper lower semicontinuous function 
$f:\R^n\rightarrow\R\cup\{+\infty\}$ has the Kurdyka-\L ojasiewicz property at 
any noncritical point $\bar{x}\in \R^n$. Indeed, Lemma \ref{semicontinuity} 
yields the existence of $c>0$ such that
$$\dist(0,\partial f(x))\geq c>0$$
whenever $x\in B(\bar{x},c/2)\cap[f(\bar x)-c/2<f<f(\bar x)+c/2]$; 
$\varphi(s)=c^{-1}s$ is then a suitable concave function.\\
(c) When a function has the Kurdyka-\L ojasiewicz property, it is important to 
have estimations of $\eta,U,\varphi$. We shall see for instance that many
convex functions satisfy the above property with $U=\R^n$ and $\eta=+\infty$. 
The determination of tight bounds for the nonconvex case is a lot more 
involved.}
\end{remark}

\subsection{Convergence to a critical point and other convergence results}

This section is devoted to the convergence analysis of the proximal 
minimization algorithm introduced in Section \ref{algo}. It provides the main 
mathematical results of this paper. Applications and examples are given in 
Section \ref{S.Appli}.  

Related previous work may be found in \cite{abrs,ARS}; but there, the setting 
is convex with a quadratic coupling $Q$, and the mathematical analysis relies 
on the monotonicity of the convex subdifferential operators. In~\cite{attbol} {\L}ojasiewicz
 inequality is used to derive the convergence of the usual proximal method but 
alternating algorithms are not considered.   

Let us introduce the notations: $z_k=(x_k,y_k)$, $l_k=L(z_k)$, 
$\bar z=(\bar x,\bar y)$, $\bar l=L(\bar z)$.

\begin{theorem}[pre-convergence result]\label{L.conv}
Assume that $L$ satisfies $(\cal H)$, $(\mathcal H_1)$ and has the 
Kurdyka-\L oja\-siewicz property at $\bar z=(\bar x,\bar y)$. Denote by $U$, 
$\eta$ and $\varphi:[0,\eta):\rightarrow\R$ the objects appearing in 
(\ref{Loja}), relative to $L$ and $\bar z$. Let $\rho>0$ be such that 
$B(\bar z,\rho)\subset U$. 
  
Let $(z_k)$ be a sequence generated by the alternating algorithm  (\ref{prox1}, 
\ref{prox2}), with $z_0$ as an initial point. Let us assume  that  
\begin{equation}\label{G}
\bar l<l_k<\bar l+\eta,\;\forall k\geq 0,
\end{equation}

\noindent and
\begin{equation}\label{H}
M\varphi(l_0-\bar l)+2\sqrt{2r_+}\sqrt{l_0-\bar l}+\|z_0-\bar z\|<\rho
\end{equation}
\noindent with $M=2r_+(C+1/r_-)$ where $C$ is a Lipschitz constant for $\nabla Q$ on 
$B(\bar z,\sqrt{2}\rho)$.

Then, the sequence $(z_k)$ converges to a critical point of $L$ and
the following estimates hold:  $\forall k\geq 0$
\begin{eqnarray}
&(i)& z_k \in B(\bar{z},\rho)  \label{est1}\\
&(ii)& \sum_{i=k+1}^\infty\|z_{i+1}-z_i\|
  \leq M\varphi(l_k-\bar l)+\sqrt{2r_+}\sqrt{l_k-\bar l}.  \label{est2}
\end{eqnarray}
\end{theorem}

\noindent The meaning of the theorem is roughly the following: a sequence 
$(z_k)$ that starts in the neighborhood of a point $\bar z$ (as given in 
(\ref{H})) and that does not improve $L(\bar z)$ (as given in (\ref{G})) 
converges to a critical point near $\bar z$. 

\smallskip\noindent
{\bf Proof.} The idea of the proof is in the line of \L ojasiewicz' argument: 
it consists mainly in proving -- and then using -- that $\varphi\circ L$ is a 
Liapunov function with decreasing rate close to $||z_{k+1}-z_k||$.

\noindent With no loss of generality one may assume that $L(\bar z)=0$ 
(replace if necessary $L$ by $L-L(\bar z)$). 

\noindent With (\ref{descent}) we have for $i\geq 0$:
\begin{equation}\label{un}
l_i-l_{i+1}\geq\frac{1}{2r_+}\|z_{i+1}-z_i\|^2.
\end{equation}
\noindent But $\varphi'(l_i)$ makes sense in view of (\ref{G}), and $\varphi'(l_i)>0$; hence
$$\varphi'(l_i)(l_i-l_{i+1})\geq\frac{\varphi'(l_i)}{2r_+}\|z_{i+1}-z_i\|^2.$$
Owing to $\varphi$ being concave, we have further:
\begin{equation}\label{deux}
\varphi(l_i)-\varphi(l_{i+1})
\geq\frac{\varphi'(l_i)}{2r_+}\|z_{i+1}-z_i\|^2,\;\forall i\geq 0.
\end{equation}

\noindent Let us first check (i) for $k=0$ and $k=1$.
In view of (\ref{H}), $z_0$ lies in $B(\bar z,\rho)$ .

\noindent As to $z_1$, (\ref{un}) yields in particular
\begin{equation}\label{z1}
\frac{1}{2r_+}\|z_1-z_0\|^2\leq l_0-l_1\leq l_0.
\end{equation}
\noindent Hence:
\begin{equation}\label{z1mzb}
\|z_1-\bar z\|\leq\|z_1-z_0\|+\|z_0-\bar z\|
  \leq\sqrt{2r_+}\sqrt{l_0}+\|z_0-\bar z\|.
\end{equation}
\noindent Hence $z_1$ lies in $B(\bar z,\rho)$ in view of (\ref{H}).

\noindent Let us now prove by induction that $(x_k,y_k)\in B(\bar z,\rho)$ for all 
$k\geq0$. 

\noindent This being true for $k\in\{0,1\}$, let us assume that it holds up to some 
$k\geq 1$.

\noindent For $0\leq i\leq k$, since $z_i\in B(\bar z,\rho)$ and $0<l_i<\eta$, 
we can write the Kurdyka-Lojasiewicz inequality at $z_i$:
$$\varphi'(l_i)\,\dist(0,\partial L(z_i))\geq 1.$$
\noindent Recall, with lemma \ref{basic} (iii), that  
$$(x_i^*,y_i^*)=
  -(\nabla_x Q(x_i,y_{i-1})-\nabla_x Q(x_i,y_i),0)
  -\left(\frac{1}{\lambda_{i-1}}(x_i-x_{i-1}),
    \frac{1}{\mu_{i-1}}(y_i-y_{i-1})\right)$$
\noindent is an element of $\partial L(x_i,y_i)$.
Hence we have for $1\leq i\leq k$:
\begin{equation}\label{gnagna}
\varphi'(l_i)\|(x_i^*,y_i^*)\|\geq 1.
\end{equation}

\noindent Let us examine $\|(x_i^*,y_i^*)\|$, for $1\leq i\leq k$. On the one 
hand
$$
\left\|\left(\frac{1}{\lambda_{i-1}}(x_i-x_{i-1}),
  \frac{1}{\mu_{i-1}}(y_i-y_{i-1})\right)\right\|
  \leq\frac{1}{r_-}\|z_i-z_{i-1}\|.$$
On the other hand, let us observe
$$\|(x_i,y_{i-1})-(\bar x,\bar y)\|^2
  =\|x_i-\bar x\|^2+\|y_{i-1}-\bar y\|^2
  \leq\|z_i-\bar z\|^2+\|z_{i-1}-\bar z\|^2
  \leq 2\rho^2.$$
Hence $(x_i,y_{i-1})$, and $z_i=(x_i,y_i)$, lie in $B(\bar z,\sqrt{2}\rho)$; 
which allows us to apply the Lipschitz inequality between these points
$$\|\nabla_x Q(x_i,y_i)-\nabla_x Q(x_i,y_{i-1})\|
  \leq C\|y_i-y_{i-1}\|
  \leq C\|z_i-z_{i-1}\|.$$
Hence, for $1\leq i\leq k$
\begin{equation}\label{maj}
\|(x_i^*,y_i^*)\|\leq(C+1/r_-)\|z_i-z_{i-1}\|.
\end{equation}
Now (\ref{gnagna}) yields
$$\varphi'(l_i)\geq\frac{1}{(C+1/r_-)}\|z_i-z_{i-1}\|^{-1},
  \;1\leq i\leq k.$$
And (\ref{deux}) yields
$$\varphi(l_i)-\varphi(l_{i+1})
  \geq\frac{1}{M}\frac{\|z_{i+1}-z_i\|^2}{\|z_i-z_{i-1}\|},\;1\leq i\leq k.$$
We rewrite the previous inequality in the following way
$$\|z_i-z_{i-1}\|^{1/2}
  \left(M(\varphi(l_i)-\varphi(l_{i+1}))\right)^{1/2}
  \geq\|z_{i+1}-z_i\|.$$
Hence (recall $ab\leq(a^2+b^2)/2$)
\begin{equation}\label{ineg}
\|z_i-z_{i-1}\|+M(\varphi(l_i)-\varphi(l_{i+1}))
  \geq 2\|z_{i+1}-z_i\|.
\end{equation}
This inequality holds for $1\leq i\leq k$; let us sum over $i$
$$\|z_1-z_0\|+M(\varphi(l_1)-\varphi(l_{k+1}))
  \geq\sum_{i=1}^k\|z_{i+1}-z_i\|+\|z_{k+1}-z_k\|.$$
Hence, in view of the monotonicity properties of $\varphi$ and $l_k$
$$\|z_1-z_0\|+M\varphi(l_0)\geq\sum_{i=1}^k\|z_{i+1}-z_i\|.$$
We finally get
$$\|z_{k+1}-\bar z\|
  \leq\sum_{i=1}^k\|z_{i+1}-z_i\|+\|z_1-\bar z\|
  \leq M\varphi(l_0)+\|z_1-z_0\|+\|z_1-\bar z\|;$$
which entails $z_{k+1}\in B(\bar z,\rho)$ in view of (\ref{z1}, \ref{z1mzb}) 
and (\ref{H}). That completes the proof of (i).

\noindent Indeed inequality (\ref{ineg}) holds for $i\geq 1$; let us sum it for $i$ 
running from some $k$ to some $K>k$
$$\|z_k-z_{k-1}\|+M(\varphi(l_k)-\varphi(l_{K+1}))
  \geq\sum_{i=k}^K\|z_{i+1}-z_i\|+\|z_{K+1}-z_K\|.$$
Hence
$$\|z_k-z_{k-1}\|+M\varphi(l_k)
  \geq\sum_{i=k}^K\|z_{i+1}-z_i\|.$$
Letting $K\rightarrow\infty$ yields
\begin{equation}\label{preconv}
\sum_{i=k}^\infty\|z_{i+1}-z_i\|\leq M\varphi(l_k)+\|z_k-z_{k-1}\|.
\end{equation}
We conclude with (\ref{un}), and that will prove point (ii):
$$\sum_{i=k}^\infty\|z_{i+1}-z_i\|
  \leq M\varphi(l_k)+\sqrt{2r_+}\sqrt{l_{k-1}}
  \leq M\varphi(l_{k-1})+\sqrt{2r_+}\sqrt{l_{k-1}}.$$

\noindent This clearly implies that $(z_k)$ is a convergent sequence.  As a 
consequence of Proposition \ref{GEN}, we obtain that its limit is a critical point of $L$.
$\hfill\Box$

\smallskip
This theorem has two important consequences

\begin{theorem}[convergence]\label{T.bounded} Assume that $L$ satisfies 
$(\mathcal H)$, $(\mathcal H_1)$ and has the Kurdyka-\L ojasiewicz property at 
each point of the domain of $f$. 

Then:

-- either  $||(x_k,y_k)||$ tends to infinity 

-- or $(x_k-x_{k-1},y_k-y_{k-1})$ is $l^1$, {\it i.e.}
$$\sum_{k=1}^{+\infty} \|x_{k+1}-x_k\|+\|y_{k+1}-y_k\|<+\infty,$$ 
and, as a consequence, $(x_k,y_k)$ converges to a critical point of $L$.
\end{theorem}
{\bf Proof.} Assume that $\|(x_k,y_k)\|$ does not tend to infinity and let 
$\bar{z}$ be a limit-point of $(x_k,y_k)$ for which we denote 
by $\rho,\eta,\varphi$ the associated objects as defined in (\ref{Loja}). 
Note that  Proposition~\ref{GEN} implies that $\bar{z}$ is critical and that 
$l_k=L(x_k,y_k)$  converges to $L(\bar{z})$. 

 If there exists an integer 
 $k_0$ for which $L(x_{k_0},y_{k_0})=L(\bar{z})$, it is straightforward to check (recall (\ref{descent})) that $(x_k,y_k)=(x_{k_0},y_{k_0})$ for all $k\geq k_0$, so that
 $(x_{k_0},y_{k_0})=\bar{z}$. 
 We may thus assume that $L(x_k,y_k)>L(\bar{z})$.

 Since 
$\max(\varphi(l_k-L(\bar{z})),\|z_k-\bar{z}\|)$ admits $0$ as a cluster point, 
we obtain the existence of $k_0\geq 0$ such that (\ref{H}) is fulfilled with $z_{k_0}$ as a new initial point. The 
conclusion is then a consequence of Theorem~\ref{L.conv}.~$\hfill\Box$
\begin{remark}{\rm (a) Several standard assumptions automatically guarantee
 the boundedness of the sequence $(x_k,y_k)$, hence its convergence: 

-- {\em One-sided coercivity implies convergence.} Assuming that $f$ (or $g$) 
has compact lower-level sets and that $Q(x,y)=\frac{1}{2}||x-y||^2$, implies
that the sequence $(x_k,y_k)$ is bounded (use Lemma~\ref{basic}).

-- {\em Convexity implies convergence.} Assume that $f,g$ are convex and that 
$Q$ is of the form $Q(x,y)=||Ax-By||$ where $A:\R^m\rightarrow \R^p$ and
$B:\R^n\rightarrow \R^p$ are linear mappings. If $L$ has at least a minimizer, 
then the sequence $(x_k,y_k)$ is bounded (see \cite{abrs}).\\
(b) Theorem \ref{T.bounded} gives new insights into convex alternating 
methods: first it shows that the finite length property is satisfied by many 
convex functions ({\it e.g.} convex definable functions, see next sections), 
but it also relaxes the quadraticity assumption on $Q$ that is required by 
the alternating minimization algorithm of \cite{abrs,ARS}.} 
\end{remark}

\medskip

The following result should not be considered as a classical result of local
convergence: we do not assume here any standard nondegeneracy conditions like 
for instance uniqueness of the minimizers, second-order conditions or 
transversality conditions. 

\begin{theorem}[local convergence to global minima]\label{T.local}
Assume that $L$ satisfies $(\cal H)$, $(\mathcal H_1)$ and has the 
Kurdyka-\L ojasiewicz property at $(\bar x,\bar y)$, a global minimum point of 
$L$. Then there exist $\epsilon$, $\eta$ such that 
$$\|(x_0,y_0)-(\bar x,\bar y)\|<\epsilon,\;\;
  \min L<L(x_0,y_0)<\min L+\eta$$
implies that the sequence $(x_k,y_k)$ starting from $(x_0,y_0)$ has the finite 
length property and converges to $(x^*,y^*)$ with 
$L(x^*,y^*)=\min L$.
\end{theorem}
{\bf Proof.} A straightforward application of Theorem \ref{L.conv} yields the 
convergence of $(x_k,y_k)$ to some $(x^*,y^*)$, a critical point of $L$ with 
$L(x^*,y^*)\in[\min L,\min L+\eta)$. Now, if $L(x^*,y^*)$ were not equal to 
$L(\bar x,\bar y)$ then the Kurdyka-Lojasiewicz inequality would entail 
$\varphi'(L(x^*,y^*)-L(\bar x,\bar y))\,\dist(0,\partial L(x^*,y^*))\geq 1$, a 
clear contradiction since $0\in\partial L(x^*,y^*)$.$\hfill\Box$  

\bigskip

The convergence rate result that follows is motivated by the analysis of 
problems involving semialgebraic or subanalytic data. These are very common 
problems (see Section \ref{S.Appli}).

\begin{theorem}[rate of convergence]\label{T.rate}
Assume that $L$ satisfies $(\cal H)$, $(\mathcal H_1)$. Assume further that 
$(x_k,y_k)$ converges to $(x_\infty,y_\infty)$ and that $L$ has 
the Kurdyka-\L ojasiewicz property at $(x_\infty,y_\infty)$ with 
$\varphi(s)=cs^{1-\theta}$, $\theta \in [0,1),\;c>0$. Then the following 
estimations hold

(i) If $\theta=0$ then the sequence  $(x_k,y_k)_{k\in \mathbb{N}}$ converges in a 
finite number of steps.

(ii) If $\theta\in(0,\frac{1}{2}]$ then there exist $c>0$ and $\tau\in[0,1)$ 
such that
$$\|(x_k,y_k)-(x_\infty,y_\infty)\|\leq c\:\tau^k.$$

(iii) If $\theta \in (\frac{1}{2},1)$ then there exists $c>0$ such that
$$\|(x_k,y_k)-(x_\infty,y_\infty)\|\leq c\:k^{-\frac{1-\theta}{2\theta-1}}.$$
\end{theorem}

\noindent
{\bf Proof.} The notations are those of Theorem~\ref{T.local}  and for simplicity 
we assume that $l_k\rightarrow 0$. Then $L(x_\infty,y_\infty)=0$ (prop. 
\ref{GEN}(iii)). 

(i) Assume first $\theta=0$. If $(l_k)$ is stationary, then so is $(x_k,y_k)$ 
in view of lemma \ref{basic}(i). If $(l_k)$ is not stationary, then the 
Kurdyka-Lojasiewicz inequality yields for any $k$ sufficiently large \
$c\,\dist(0,\partial L(x_k,y_k))\geq 1$, a contradiction in view of lemma 
\ref{basic}(iii).  

(ii,iii) Assume $\theta>0$. For any $k\geq0$, set 
$\Delta_k=\sum_{i=k}^\infty\sqrt{\|x_{i+1}-x_i\|^2+\|y_{i+1}-y_i\|^2}$ which 
is finite by Theorem \ref{T.bounded}. Since 
$\Delta_k\geq\sqrt{\|\xk-x_\infty\|^2+\|y_k-y_\infty\|^2}$, it is sufficient
to estimate $\Delta_k$. The following is a rewriting of (\ref{preconv}) with 
these notations 
\begin{equation}\label{deltak}
\Delta_k\leq M\varphi(l_k)+(\Delta_{k-1}-\Delta_k).
\end{equation}
The Kurdyka-\L ojasiewicz inequality successively yields 
\begin{eqnarray*}
& \varphi'(l_k)\dist(0,\partial L(x_k,y_k))
= c(1-\theta)l_k^{-\theta}\dist(0,\partial L(x_k,y_k)) \geq 1 & \\
& l_k^\theta\leq c(1-\theta)\dist(0,\partial L(x_k,y_k)). &
\end{eqnarray*}
But with (\ref{maj}) we have 
$$\dist(0,\partial L(x_k,y_k))\leq\|(x_k^*,y_k^*)\|
  \leq(C+1/r_-)(\Delta_{k-1}-\Delta_k).$$
Combining the previous two inequalities we obtain for some positive $K$ 
$$\varphi(l_k)=cl_k^{1-\theta}
  \leq K(\Delta_{k-1}-\Delta_k)^\frac{1-\theta}{\theta}.$$
Finally (\ref{deltak}) gives 
$$\Delta_k\leq 
MK(\Delta_{k-1}-\Delta_k)^{\frac{1-\theta}{\theta}}+(\Delta_{k-1}-\Delta_k).$$
Sequences satisfying such inequalities  have been studied in 
\cite[Theorem 2]{attbol}. Items (ii) and (iii) follow from these results. 
$\hfill\Box$

\subsection{Convergence of alternating projection methods}\label{vonN}
In this section we consider the special but important case of bifunctions of 
the type
\begin{equation}\label{LCD}
L_{C,D}(x,y)=\delta_C(x)+\frac{1}{2}\|x-y\|^2+\delta_D(y),\;(x,y)\in\R^n,
\end{equation}
where $C,D$ are two nonempty closed subsets of $\R^n$. Notice that $L$ 
satisfies $(\mathcal H)$ and $(\mathcal H_1)$ (see Section \ref{algo}) for any 
$y_0\in D$. In this specific setting, the proximal minimization algorithm 
(\ref{prox1}, \ref{prox2}) reads
\begin{eqnarray*}
x_{k+1} & \in & \argmin\{\frac{1}{2}\|u-y_k\|^2+\frac{1}{2\lambda_k}\|u-x_k\|^2 :u\in C\}\\
y_{k+1} & \in & \argmin \{\frac{1}{2}\|v-x_{k+1}\|^2+\frac{1}{2\mu_k}\|v-y_k\|^2: v\in D\}.
\end{eqnarray*}
Thus we obtain the following {\em alternating projection algorithm}
\begin{equation*}
\left\{\begin{aligned}
x_{k+1} & \in  P_C\left(\frac{\lambda_k^{-1}x_k+y_k}{\lambda_k^{-1}+1}\right)\\
y_{k+1} & \in  P_D\left(\frac{\mu_k^{-1}y_k+x_{k+1}}{\mu_k^{-1}+1}\right).
\end{aligned}\right.
\end{equation*}
 
The following result illustrates the interest of the above algorithm for feasibility problems.

\begin{corollary}[convergence of sequences] 
Assume that the bifunction $L_{C,D}$ has the Kur\-dyka-\L oja\-sie\-wicz 
property 
at each point. Then either $\|(x_k,y_k)\|\rightarrow\infty$ as 
$k\rightarrow\infty$, or $(x_k,y_k)$ converges to a critical point of $L$.\\
\noindent
{\bf(local convergence)} Assume that the bifunction $L_{C,D}$ has the 
Kurdyka-\L ojasiewicz property at $(x^*,y^*)$ and that
$\|x^*-y^*\|=\min\{\|x-y\|:x\in C,\,y\in D\}$. If $(x_0,y_0)$ is sufficiently 
close to $(x^*,y^*)$ then the whole sequence converges to a point 
$(x_\infty,y_\infty)$ such that 
$\|x_\infty-y_\infty\|=\min\{\|x-y\|:x\in C,\,y\in D\}$.
\end{corollary}
\noindent
{\bf Proof.} The first point is due to Theorem \ref{T.bounded} while the 
second follows from a specialization of Theorem~\ref{T.local} to $L_{C,D}$.
$\hfill\Box$

\smallskip
Observe that for $\lambda_k$ and $\mu_k$ large, the method is very close to 
the von Neuman alternating projection method. Section~\ref{S:feasibility}  
lists a wide class of sets for which convergence of the sequence is ensured. 
%%%%%%%

%%% Figure 1

%%%%%%

\section{Examples and applications}\label{S.Appli}

A first simple theoretical result which deserves to be mentioned is that a 
generic smooth function satisfies the \L ojasiewicz inequality. Recall that a 
$C^2$ function $f:\R^n\rightarrow\R$ is a {\em Morse function} if for each 
critical point $\bar{x}$ of $f$, the hessian $\nabla ^2f(\bar{x})$ of $f$ at 
$\bar x$ is a nondegenerate endomorphism of $\R^n$. Morse functions can be 
shown to be generic in the Baire sense in the space of $C^2$ functions 
(see \cite{AE}).

Let $\bar{x}$ be a critical point of $f$, a Morse function.  Using the Taylor formula for $f$
and $\nabla f$, we obtain the existence of a neighborhood $U$ of $\bar x$ and 
positive constants $c_1,c_2$ for which
$$|f(x)-f(\bar{x})|\leq c_1||x-\bar{x}||^2,\; 
  ||\nabla f(x)||\geq c_2||x-\bar{x}||,$$
whenever $x\in U$. It is then straightforward to see that $f$ complies with 
$(\ref{Loja})$ with a function $\varphi$ of the form $\varphi(s)=c\sqrt{s}$, 
where $c$ is a positive constant.

\medskip
Let us now come to more practical aspects and to several concrete 
illustrations of the Kurdyka-Lojasiewicz inequality (\ref{Loja}).

\subsection{Convex examples}

A general smooth convex function may not satisfy the Kurdyka-\L ojasiewicz 
(see \cite{BDLM} for a counterexample), however, in many practical cases 
convex functions do satisfy this inequality.

\smallskip\noindent
{\bf Growth condition for convex functions}\label{ss1}: 
Consider a convex function $f$ satisfying the following growth condition:
$\exists U\mbox{ neighborhood of }\bar x$, $\eta>0$, $c>0$, $r\geq1$ such that
$$\forall x\in U\cap[\min f<f<\min f+\eta],\:
  f(x)\geq f(\bar x)+cd(x,\argmin f)^r,$$
where $\bar x\in\argmin f\neq\emptyset$. Then $f$ complies with (\ref{Loja}) 
at point $\bar x$ (for $\varphi(s)=r\;c^{-\frac{1}{r}}\;s^{\frac{1}{r}}$) on 
$U\cap[\min f<f<\min f+\eta]$ (see \cite{BDL04}).

\smallskip\noindent
{\bf Uniform convexity}:
If $f$ is uniformly convex i.e., satisfies
$$f(y)\geq f(x)+\langle x^*,y-x\rangle +K\|y-x\|^p,\;p\geq 1 $$
for all $x,y\in \R^n$, $x^*\in \partial f(x)$ then $f$ satisfies the 
Kurdyka-\L ojasiewicz inequality on $\dom f$ for 
$\varphi(s)=p\;K^{-\frac{1}{p}}\;s^{\frac{1}{p}}$.

\noindent
{\bf Proof.} Since $f$ is coercive and strictly convex, 
$\argmin f=\{\bar x\}\neq\emptyset$. Take $y\in\R^n$. By applying the uniform 
convexity property at the minimum point $\bar x$, we obtain that  
$$f(y)\geq\min f +K\|y-\bar x\|^p.$$
The conclusion follows from the preceding paragraph.

\smallskip \noindent
{\bf Tame convex functions}: Another important class of functions which is 
often met in practice is the class of convex functions which are definable in 
an o-minimal structure (e.g. semialgebraic or globally subanalytic convex 
functions). In a tame setting the Kurdyka-\L ojasiewicz inequality does not 
involve any convexity assumptions, the reader is thus referred to section 
\ref{S:omin} for an insight into this kind of results. 

\medskip
Among the huge literature devoted to convex problems, one may consult 
\cite{BB,Comb,BC,CombWajs} and their references.

\subsection{Metrically regular equations} 
Let the mapping $F:\R^n\rightarrow\R^m$ be metrically regular at some point 
$\bar x\in\R^n$, namely (cf. \cite{Rock98,Ioffe}) there exist a neighborhood 
$V$ of $\bar x$ in $\R^n$, a neighborhood $W$ of $F(\bar x)$ in $\R^m$ and a 
positive $k$ such that 
$$x\in V,\;y\in W\;\Rightarrow\;\dist(x,F^{-1}(y))\leq k\;\dist(y,F(x)).$$
The coefficient $k$ is a measure of stability under perturbation of equations 
of the type $F(\bar{x})=\bar{y}$.

Let $C\subseteq\R^m$ and consider the problem of finding a point $x\in\R^n$ 
satisfying $F(x)\in C$, so that we are concerned with metric regularity in 
constraint systems.

For the sake of simplicity, we assume that $F$ is $C^1$ on a neighborhood of 
$\bar x$ and that $C$ is convex, closed, nonvoid.

Solving our constraint system amounts to solving  
$$\min\left\{ f(x):=\frac{1}{2}\dist^2(F(x),C):x\in\R^n\right\},$$
so that $f$ is now the object under study. The function $f$ is differentiable 
around $\bar x$ with gradient 
$$\nabla f(x)=D^*F(x)[F(x)-P_C(F(x))],$$ 
where $D^*F(x)$ denotes the adjoint of the Fr\'echet derivative of $F$ at $x$. 

The metric regularity of $F$ at $\bar x$ entails that 
$D^*F(\bar x):\R^m\mapsto\R^n$ is one-to-one, with $\|[D^*F(x)]^{-1}\|\leq k$ 
(\cite[Th. 9.43]{Rock98}\cite[Ch. 3, Th. 3]{Ioffe}). Fix some $k'>k$. In view 
of the continuity of $\|[D^*F(x)]^{-1}\|$ around $\bar x$, there exists a 
neighborhood $U$ of $\bar x$ such that: 
$x\in U\;\Rightarrow\;\|[D^*F(x)]^{-1}\|\leq k'$. We then deduce from the 
expression above for $\nabla f(x)$  
$$x\in U\;\Rightarrow\;\|\nabla f(x))\|\geq\frac{1}{k'}\|F(x)-P_C(F(x))\|
  =\frac{1}{k'}\dist(F(x),C)=\frac{1}{k'}\sqrt{2f(x)}.$$
Hence
$$x\in U,\;f(\bar x)<f(x)\;\Rightarrow
  \|\nabla f(x))\|\geq\frac{1}{k'}\sqrt{2(f(x)-f(\bar x))}.$$
In other words, $f$ satisfies the Kurdyka-\L ojasiewicz inequality at all 
point $\bar x$ where $F$ is metrically regular, with 
$\varphi:s\in[0,+\infty[\mapsto k'\sqrt{2s}.$

It is important to observe here that $\argmin f$ can be a {\em continuum} (in 
that case $f$ is not a Morse function). This can easily be seen by taking 
$m<n$, $b\in\R^m$, a full rank matrix $A\in\R^{m\times n}$ and $F(x)=Ax-b$ for 
all $x$ in $\R^n$.

For links between metric regularity and the Kurdyka-Lojasiewicz inequality 
see \cite{BDLM}.

\subsection{Tame functions}\label{S:omin}

As it was  emphasized in the introduction, tame sets and functions provide
 a vast field of applications of our main results.

\bigskip

\noindent
{\bf Semialgebraic functions }: Recall that a subset of $\R^n$ is called semialgebraic if it can be written as a finite
 union of sets of the form
$$\{x\in \R^n: p_{i}(x)=0,\; q_{i}(x)<0,\;i=1,\ldots,p\},$$ 
where $p_{i},q_{i}$ are real polynomial functions.

A function $f:\R^n\rightarrow \R\cup\{+\infty\}$ is semialgebraic if its graph 
is a semialgebraic subset of $\R^{n+1}$. Such a function satisfies the 
Kurdyka-\L ojasiewicz property (see \cite{BDL04,BDLS}) with 
$\varphi(s)=cs^{1-\theta}$, for some $\theta\in[0,1)\cap\mathbb{Q}$ and some 
$c>0$. This nonsmooth result generalizes the famous 
\L ojasiewicz inequality for real-analytic functions \cite{Loja63}.
Stability properties of semialgebraic functions are numerous (see e.g. 
\cite{BR1990,BcR1998}), the following few facts might help the reader to 
understand how they impact Optimization matters:

-- finite sums and products of semialgebraic functions are semialgebraic;

-- scalar products are semialgebraic;

-- indicator functions of semialgebraic sets are semialgebraic;

-- generalized inverse of semialgebraic mappings are semialgebraic;

-- composition of semialgebraic functions or mappings are semialgebraic;

-- functions of the type $\R^n \ni x\rightarrow f(x)=\sup_{y\in C} g(x,y)$ 
(resp. $\R^n \ni x\rightarrow f(x)=\inf_{y\in C} g(x,y)$) where $g$ and $C$ 
are semialgebraic are semialgebraic.

Matrix theory provides a long list of semialgebraic objects (see \cite{mallew}): positive 
semidefinite matrices, Stiefel manifolds (spheres, orthogonal group; see e.g. 
\cite{Edel}), constant rank matrices... Let us provide now some concrete 
problems for which our algorithm can be implemented effectively.
 
\begin{example}[Rank reduction of correlation matrices]
{\rm The following is a standard problem: being given a symmetric matrix $A$ and an 
integer $d\in\{1,\ldots,n-1\}$, we wish to find a correlation matrix $Z$ of 
small rank and as close as possible to $A$. An abstract formulation of the 
problem leads to the following
$$\min \{F(Z-A) : Z\in\s_n^+,\,\diag Z=I_n,\,\rang Z\leq d\},$$
where $F:\R^{n\times n}\rightarrow \R$ is a smooth semialgebraic function ({\it e.g.} a seminorm),
 $\s_n^+$ is the cone of positive semidefinite $n$ matrices and $I_n$ is the $n\times n$ identity matrix. 
The reformulation of the rank reduction problem for correlation matrices we 
use here, is due to Grubi\v si\'c-Pietersz \cite{GrubP} who transform the 
original problem into the minimization of a functional on a matrix manifold 
(for more on optimization on matrix manifolds see \cite{AMS}) . 
%Let 
%$G:Y\in\R^{n\times d}\rightarrow G(Y)\in\R$ be a smooth function defined on 
%the space of $n\times d$ real matrices ($0<d<n$), identified with 
%$\R^{n\times d}$. 
Let $Y_i$, $i\in\{1,\ldots,n\}$, denote line $i$ of 
$Y\in\R^{n\times d}$. Consider the following manifold of $\R^{n\times d}$ 
$$\mathcal C\,=\,\left\{Y\in\R^{n\times d}:\,
  Y_i\in S^{i-1}\times\{0\}^{d-i},\,\forall i\in\{1,\ldots,d\};\,
  Y_i\in S^{d-1},\,\forall i\in\{d+1,\ldots,n\}\right\},$$
where $S^0=\{1\}$ and $S^{i-1}$ denotes the unit sphere of $\R^i$ for 
$i\in\{2,\ldots,d\}$. The set $\cal C$ is a real-algebraic manifold, and, due 
to its simple structure (it is a product of low dimensional spheres), the 
projection operator onto $\cal C$ is simply given as a product of projections 
on spheres. Denote by $Y^T$ the transpose matrix of $Y$. Since   $\{Z\in \s_n^+: \diag Z=I_n, \rang Z\leq d\}=\{YY^T: Y\in {\cal C}\}$ (see \cite{GrubP}),   
the correlation matrix rank reduction problem can be reformulated as
\begin{equation}\label{rang}
\min\;\{\,G(Y)\,:\:Y\in\mathcal C\,\},
\end{equation}
where $G(Y)=F(YY^T-A)$ . 
In \cite{GrubP} this problem is addressed by Newton's method or a conjugate 
gradient algorithm on the manifold $\cal C$. Newton's method enjoys its usual 
local quadratic convergence, while theoretical results about the conjugate 
gradient are sparse in spite of good numerical evidence.

We propose the following relaxation of problem (\ref{rang}) 
\begin{equation}\label{relax}
\min\;\left\{\,\delta_{\mathcal C}(X)
  +\frac{\rho}{2} 
    \sum_{\stackrel{\scriptstyle i=1,\ldots,n}{j=1,\ldots,d}}
      (X_{ij}-Y_{ij})^2
  +G(Y)\,:\:X,Y\in\R^{n\times d}\right\},
\end{equation}
where $\rho>0$ is a penalization parameter. Since the 
functional $L(X,Y)=\delta_{\mathcal C}(X)+\sum_{i,j}(X_{ij}-Y_{ij})^2+G(Y)$ is 
semialgebraic, the problem (\ref{relax}) falls into the context of our 
study. Algorithm (\ref{algo0}) applied to (\ref{relax}) generates a sequence 
$(X_k,Y_k)$ that converges ($\cal C$ is bounded) to a critical point 
$(\bar X,\bar Y)$ of $L$ with $\bar X\in\mathcal C$. The algorithm reads
\begin{equation}\label{algogo}
\left\{
\begin{aligned}
 X_{k+1} & \in\mbox{argmin}
  \left\{\frac{\rho}{2}
     \sum_{\stackrel{\scriptstyle i=1,\ldots,n}{j=1,\ldots,d}}
      (X_{ij}-Y_{k,ij})^2+
    \frac{1}{2\lambda_k}
     \sum_{\stackrel{\scriptstyle i=1,\ldots,n}{j=1,\ldots,d}}
      (X_{ij}-X_{k,ij})^2\ :\ X\in C \right\}\\
 Y_{k+1} & \in \mbox{argmin}
  \left\{G(Y)+
    \frac{\rho}{2}
     \sum_{\stackrel{\scriptstyle i=1,\ldots,n}{j=1,\ldots,d}}
      (Y_{ij}-X_{k+1,ij})^2+
    \frac{1}{2\mu_k}
     \sum_{\stackrel{\scriptstyle i=1,\ldots,n}{j=1,\ldots,d}}
      (Y_{ij}-Y_{k,ij})^2\ :\ Y\in\R^{n\times d}\right\}
\end{aligned}\right.
\end{equation}
The second minimization is a proximal step for the function 
$Y\rightarrow\frac{\rho}{2}\sum_{i,j}(X_{k+1,ij}-Y_{ij})^2+G(Y)$ (which may 
be convex, for large $\rho$). The first minimization reduces to $n$ uncoupled 
minimizations in $\R^d$ which give the $n$ lines $X_{k+1,i}$ of matrix 
$X_{k+1}$:
$$
X_{k+1,i} \mbox{ is a  minimizer of }
 \R^d \ni X_i\rightarrow  \sum_{j=1,\ldots,d}\left(X_{ij}-
    \frac{\rho\lambda_k Y_{k,ij}+X_{k,ij}}{\rho\lambda_k +1}\right)^2
  ,
$$
where $X_i$ is subject to the constraints: $X_i\in S^{i-1}\times\{0\}^{d-i}$ if 
$i\in\{1,\ldots,d\}$, and $X_i\in S^{d-1}$ if $i\in\{d+1,\ldots,n\}$.\\

 Each 
solution line $X_{k+1,i}$ is simply the projection of line  
$\rho\lambda_k Y_{k,i}+X_{k,i}$ of matrix $\rho\lambda_k Y_k+X_k$ onto the 
unit sphere of $\R^d$ if $d\leq i\leq n$, and onto a unit sphere of lower 
dimension if $2\leq i\leq d-1$. The second step is therefore solved in closed form.}
\end{example}

\noindent
{\bf Functions definable in an o-minimal structure over $\R$.} Introduced in 
\cite{Dries} these structures can be seen as an axiomatization of the  
qualitative properties of semialgebraic sets. 

\begin{definition}\label{D:omin}
Let ${\cal O}=\{{\cal O}_n\}_{n\in \N}$ be  such that each
 ${\cal O}_n$ is a collection of subsets of $\R^n$. The family ${\cal O}$ is an o-minimal structure over $\R$, if it satisfies
 the following axioms:
\begin{itemize}
\item[(i)] Each ${\cal O}_n$ is a boolean algebra. Namely 
$\emptyset\in{\cal O}_n$ and for each $A,B$ in ${\cal O}_n$, $A\cup B$, 
$A\cap B$ and $\R^n\setminus A$ belong to ${\cal O}_n$.
\item[(ii)] For all $A$ in ${\cal O}_n$, $A\times \R$ and $\R \times A$ belong to  ${\cal O}_{n+1}$.
\item[(iii)] For all $A$ in ${\cal O}_{n+1}$, $\Pi(A):=\{(x_1,\ldots,x_n)\in \R^n: (x_1,\ldots,x_n,x_{n+1})\in A\}$ belongs to ${\cal O}_n$.
\item[(iv)] For all $i\neq j$ in $\{1,\ldots,n\}$, $\{(x_1,\ldots,x_n)\in \R^n:x_i=x_j\}\in {\cal O}_n$.
\item[(v)] The set $\{(x_1,x_2)\in \R^2:\;x_1<x_2\}$ belongs to ${\cal O}_2$.
\item[(vi)] The elements of ${\cal O}_1$ are exactly finite unions of intervals.
\end{itemize}
\end{definition}

Let  $\cal O$ be an o-minimal structure. A set $A$ is said to be definable (in 
$\cal O$), if $A$ belongs to $\cal O$. A point to set mapping 
$F:\R^n\rightrightarrows\R^m$ (resp. a real-extended-valued function 
$f:\R^n\rightarrow\R\cup\{+\infty\}$) is said to be definable if its graph is 
a definable subset of $\R^n\times\R^m$ (resp. $\R^n\times\R$).

As announced in the introduction, definable sets and mappings share many 
properties of semialgebraic objects. Let $\cal O$ be an o-minimal structure. 
Then,\\
- finite sums of definable functions are definable;\\
- indicator functions of definable sets are definable;\\
- generalized inverses of definable mappings are definable;\\
- compositions of definable functions or mappings are definable;\\
- functions of the type $\R^n\ni x\rightarrow f(x)=\sup_{y\in C}g(x,y)$ 
(resp. $\R^n\ni x\rightarrow f(x)=\inf_{y\in C}g(x,y)$) where $g$ and $C$ are 
definable, are definable.

Due to their dramatic impact on several domains in mathematics, these 
structures are being intensively studied. One of the interests of such 
structures in optimization is due to the following nonsmooth extension of the 
Kurdyka-\L ojasiewicz inequality.

\begin{theorem}[\cite{BDLS}]\label{bdls} 
{\em Any proper lower semicontinuous function 
$f:\R^n\rightarrow \R\cup\{+\infty\}$ which is definable in an o-minimal 
structure $\cal O$ has the Kurdyka-\L ojasiewicz property at each point of 
$\dom\partial f$. Moreover the function $\varphi$ appearing in (\ref{Loja}) is 
definable in $\cal O$.}
\end{theorem}

The concavity of the function $\varphi$ is not stated explicitly in \cite{BDLS}. 
The proof of that fact is however elementary: it relies on the following  
fundamental result known as the monotonicity Lemma (see \cite{Dries}). 

(Monotonicity Lemma) {\em Let $k$ be an integer and 
$f:I\subseteq\R\rightarrow \R$ a function definable in some o-minimal 
structure. Then there exists a finite partition of $I$ into intervals 
$I_1,\ldots,I_p$ such that the restriction of $f$ to each $I_i$ is $C^k$ and 
either strictly monotone or constant.}

\noindent
{\bf Proof of Theorem \ref{bdls} [concavity of $\varphi$]} Let $x$ be a critical 
point of $f$ and $\varphi$ a \L ojasiewicz function of $f$ at $x$ (see 
(\ref{Loja})). As recalled above, such a function $\varphi$ exists and is 
definable in $\cal O$ moreover. When applied to $\varphi$, the monotonicity lemma 
yields the existence of $r>0$ such that $\varphi$ is $C^2$ with either 
$\varphi''\leq 0$ or $\varphi''>0$ on $(0,r)$ (apply the monotonicity lemma to 
$\varphi''$). If $\varphi''>0$ then $\varphi'$ is increasing. The function $\varphi$ in 
inequality (\ref{Loja}) can be replaced by $\psi(s)=\varphi'(r)s$ for all 
$s\in (0,r)$, since $\psi'(s)=\varphi'(r)>\varphi'(s)$.$\hfill\Box$.

\medskip

The following examples of o-minimal structures illustrate the considerable 
wealth of such a concept.

\smallskip

\noindent
{\bf Semilinear sets}
A subset of $\R^n$ is called semilinear if it is a finite union of sets of the form
$$\{x\in \R^n: \langle a_{i},x\rangle=\alpha_{i},\; \langle b_{i},x\rangle <\beta_{i},\;i=1,\ldots,p\},$$
where $a_{i}, b_{i}\in \R^n$ and $\alpha_{i}, \beta_{i}\in \R$.
 One can easily establish that such a structure is an o-minimal structure.

Besides the function $\varphi$ is of the form $\varphi(s)=cs$ with $c>0$.

\smallskip

\noindent
{\bf Real semialgebraic sets} By Tarski quantifier elimination theorem
 the class of semialgebraic sets is an o-minimal structure \cite{BR1990,BcR1998}.

\smallskip

\noindent
{\bf Globally subanalytic sets }(Gabrielov \cite{gabrielov}) There exists an o-minimal
structure, denoted by~$\R_{\rm{an}}$, that contains all sets of
the form $\{(x,t)\in [-1,1]^n\times\R: f(x)=t\}$ where
$f:[-1,1]^n\rightarrow \R$ ($n\in \N$) is an analytic function that can be
extended analytically on a neighborhood of the square~$[-1,1]^n$.
The sets belonging to this structure are called {\em globally
subanalytic sets}. As for the semialgebraic class, the function appearing in Kurdyka-\L ojasiewicz inequality is of the form $\varphi(s)=s^{1-\theta}$,
 $\theta\in[0,1)\cap \mathbb{Q}$.

Let us give some concrete examples. Consider  a finite collection of real-analytic functions $f_i:\R^n\rightarrow \R$ where $i=1,\ldots,p$.

\noindent
- The restriction of the function $f_+=\max_i f_i$ (resp. $f_-=\min_i f_i$) to each $[-a,a]^n$ ($a>0$) is globally subanalytic.\\ 
- Take a finite collection of real-analytic functions $g_j:\R^n\rightarrow \R$ and set 
$$C=\{x\in \R^n:f_i(x)=0, \ g_j(x)\leq 0\}.$$
If $C$ is a bounded subset then it is globally subanalytic. Let now $G:\R^m\times \R^n\rightarrow \R$ be a real
 analytic function. When $C$ is bounded, the restriction of the following function to each $[-a,a]^n$ ($a>0$) is globally subanalytic:
$$f(x)=\max_{y\in C} G(x,y).$$

\smallskip

\noindent
{\bf Log--exp structure} (Wilkie, van der Dries) \cite{Wilkie,Dries} There 
exists an o-minimal structure containing $\R_{\rm{an}}$ and the graph of 
$\exp\,:\R\rightarrow\R$. This huge structure contains all the aforementioned 
structures. One of the surprising specificity of such a structure is the 
existence of ``infinitely flat" functions like 
$x\rightarrow\exp(-\frac{1}{x^2})$. 

Many optimization problems are set in such a structure. When it is possible, 
it is however important to determine the minimal structure in which a problem 
is definable. 

This can for instance have an impact on the convergence analysis and in 
particular on the knowledge of convergence rates (see Theorem \ref{T.rate}).    

\medskip

\begin{example}\label{compressed}[Compressive sensing]
{\rm A current active trend in signal recovery aims at selecting a sparse 
({\em i.e.} with many zero components) solution of an underdetermined linear 
system. So the following problem comes under consideration (see 
\cite{CWB,BDE,donoho})
\begin{equation}\label{l0}
\min\{\|x\|_0:Ax=b\}
\end{equation}
where $\|x\|_0$ is the so-called $\ell^0$ ``norm'' which counts the nonzero 
components of $x\in\R^n$, $A$ is an $m\times n$ matrix ($m<n$) and $b\in\R^m$; 
the set of contraints $\mathcal C=\{x\in\R^n:Ax=b\}$ is supposed nonvoid. Due 
to its extremely combinatorial nature, problem (\ref{l0}) is untractable and 
signal theory scientists have to turn to alternatives. We briefly present one
of them below, in connection with our work.

In \cite{CWB}, the $\ell^0$ norm is approximated by a weighted 
$\ell^1$ norm  
\begin{equation}\label{l1}   
\min\left\{\sum_{i=1}^n w_i\mid x_i\mid:Ax=b\right\}.
\end{equation}
For fixed positive weights $w_i$, problem (\ref{l1}) classically boils down to 
a linear program whose solution may well bear little relation with that of 
(\ref{l0}). So actually, in this formulation, the weights are somehow part of 
the unknowns and are to be chosen at best. The reweighted $\ell^1$ 
minimization algorithm used in \cite{CWB} updates, at step $k$, the solution 
candidate $x^k$ and the weight vector $w^k=(w^k_1,\ldots,w^k_n)$ according to 
the following rule. 
Fix $\eps>0$:\\
-- $w^k$ being known, solve (\ref{l1}) with $w^k$ in place of $w$ to obtain 
$x^{k+1}$;\\
-- update each component of $w^k$ to $w_i^{k+1}=1/(\mid x_i^{k+1}\mid+\eps)$.

\smallskip 

We first observe that the reweighted algorithm can be seen as an alternating minimization method. 
For $(x,w)\in\R^n\times\R^n$ define $f_0(x)=\delta_{\mathcal C}(x)$, 
$Q_0(x,w)=\sum_{i=1}^n(\mid x_i\mid+\eps)w_i$, and $g_0(w)=-\sum_{i=1}^n\log w_i$ if $w_i>0$ for all $i$ in $\{1,\ldots,n\}$, $g_0(w)=+\infty$ otherwise. 
 Set $L_0(x,w)=f_0(x)+Q_0(x,w)+g_0(w)$; it is not difficult to realize that the original 
reweighted $\ell^1$ minimization algorithm is {\em exactly} the alternate 
minimization of the function $L_0$. 

In order to apply a {\em proximal} alternating method to $L_0$,  the problem $\min L_0$
 can be first reformulated as a problem with a smooth coupling term.  Let us therefore split $x$ into its nonnegative and nonpositive parts: 
$x'_i=\max(x_i,0)$, $x'_{i+n}=\max(-x_i,0)$ for $i\in\{1,\ldots,n\}$. Set  
$\mathcal C'=
 \{x'\in\R^{2n} : [A,-A]x'=b,x'_i\geq 0,\forall i\in\{1,\ldots,2n\}\}$, and 
for $x'\in\R^{2n}$ define $f(x')=\delta_{\mathcal C'}(x')$, 
$Q(x',w)=\sum_{i=1}^n(x'_i+x'_{i+n}+\eps)w_i$ and 
$L(x',w)=f(x')+Q(x',w)+g_0(w)$. Then minimizing $L_0(x,w)$ is equivalent to 
minimizing $L(x',w)$, with the correspondance $x_i=x'_i-x'_{i+n}$. 

Since $L$ is definable in  the $\log-\exp$ structure, the proximal alternating minimization (\ref{algo0}) applied to $L(x',w)$ 
yields a sequence $({x'}^k,w^k)$ whose behaviour is ruled by 
Theorem~\ref{T.bounded}. 
%Due to the proximal 
%terms in the algorithm, $({x'}^k,w^k)$ is different from the sequence the 
%reweighting algorithm would produce. 
In particular $({x'}^k,w^k)$, which is 
bounded because $L$ is coercive, converges to a critical point 
$(\tilde x',\tilde w)$ of $L$ verifying
$$
\left\{
\begin{array}{ccl}
 0\in
  \partial_{x'}
   \{Q(\tilde x',\tilde w)+\delta_{\mathcal C'}(\tilde x')\}\ ,&
 \mbox{hence}&
 \tilde x'\in
  \argmin
   \{Q(x',\tilde w)+\delta_{\mathcal C'}(x')\,:\,x'\in\R^{2n}\}\\
 0=\nabla_w\{Q(\tilde x',\tilde w)+g_0(\tilde w)\}\ ,&
 \mbox{ hence }&
 \tilde w_i={\displaystyle\frac{1}{\tilde x'_i+\tilde x'_{i+n}+\eps}},\ 
  \forall i\in\{1,\ldots,n\}
 \rule[2em]{0em}{0em}\\
\end{array}\right.
$$
where the minimizing property of $\tilde x'$ stated in the first line is a 
consequence of the convexity of the function 
$x'\rightarrow Q(x',w)+\delta_{\mathcal C'}(x')$ for a fixed $w$. Now it is 
easy to see that $\tilde x'$ must be such that 
$\min(\tilde x'_i,\tilde x'_{i+n})=0$ for $i\in\{1,\ldots,n\}$. Let 
$\tilde x\in\R^n$ be defined by $\tilde x_i=\tilde x'_i-\tilde x'_{i+n}$, then 
$\tilde x$ is a solution of (\ref{l1}) with 
$w_i=\tilde w_i=1/(|\tilde x_i|+\eps)$. 
Observe also that $L$ is (globally) subanalytic on bounded boxes so that its desingularizing functions are of the form $\varphi(s)=c s^{\theta}$ with $c>0$ and $\theta\in(0,1]$. Theorem~\ref{T.rate} on the 
 rate estimation therefore applies to the sequence $(x'^k,w^k)$. Let us summarize the previous discussion:

\medskip 

\noindent
{\bf  A proximal reweighted $\ell^1$ algorithm.}$\mbox{}$\\ 
1) Fix $(x'_0,w_0)$ in 
$\R^{2n}\times \R^{m}$. Let $\lambda_k,\mu_k$ be some sequences of stepsizes 
such that $\lambda_k,\mu_k \in (r_-,r_+)$ where $r_-$, $r_+$ are positive 
bounds (with $r_-<r_+$).\\ 
2) Generate a sequence $({x'}^k,w^k)$ such that
$$\mbox{\bf (PR) }
\left\{
\begin{aligned}
 {x'}^{k+1} & =\mbox{argmin}
  \left\{\sum_{i=1}^n(x'_i+x'_{i+n}+\eps)w_i^k+
    \frac{1}{2\lambda_k}
     \sum_{i=1}^{2n}(x'_i-{x'}^k_i)^2
      \ :\ x'\in\R_+^{2n}, \, [A,-A]x'=b \right\}\\
 w^{k+1} & =\mbox{argmin}
  \left\{\sum_{i=1}^n({x'}^{k+1}_i+{x'}^{k+1}_{i+n}+\eps)w_i-
   \sum_{i=1}^n \log w_i+
    \frac{1}{2\mu_k}
     \sum_{i=1}^n(w_i-w_i^k)^2
      \ :\ w\in\R^n\right\}
\end{aligned}\right.
$$
\noindent
{\bf Convergence of (PR).} The sequence $(x'_k,w_k)$ converges to a point $(x',w)$ and there exist positive constants $d,r$ such that
$$||(x'_k,w_k)-(x',w)||\leq \frac{d}{k^r}, \; \forall k\geq1.$$
 Moreover if we set 
$x=(x'_1-x'_{1+n},\ldots,x'_i-x'_{i+n},\ldots,x'_n-x'_{2n})$, then
 $x$ is solution
 of the reweighted minimization problem (\ref{l1}) with 
$w_i=1/(|x_i|+\eps)$.\\

\smallskip

Observe finally that the first minimization in (PR) is a convex quadratic program and that the second one 
may be solved in closed form. Indeed $w^{k+1}_i$ is the positive root of the second 
degree equation
$${\displaystyle\frac{1}{\mu_k}}w^2+
 \left[({x'}^{k+1}_i+{x'}^{k+1}_{i+n}+\eps)-
  {\displaystyle\frac{w_i^k}{\mu_k}}\right]w-
 1=0.$$
The computational burden of the proximal reweighting algorithm is not much 
heavier than that of the reweighting algorithm.}
\end{example}

\subsection{Kurdyka-\L ojasiewicz inequality for feasibility problems}
Recall the function $L_{C,D}$ defined by (\ref{LCD}). Writing down the 
optimality condition we obtain 
$$\partial L_{C,D}(x,y)=\{(x-y+u,y-x+v): u\in N_C(x), v\in N_D(y)\}.$$
This implies 
\begin{equation}\label{lojaset}
\dist(0,\partial L_{C,D}(x,y))
  =\left(\dist^2(y-x,N_C(x))+\dist^2(x-y,N_D(y))\,\right)^{\frac{1}{2}}.
\end{equation}

\medskip
\subsubsection{(Strongly) regular intersection}
The following result is a reformulation of \cite[Theorem 17]{mallew} where
this fruitful concept was considered in relation with alternate projection
methods. For the sake of completeness we give a proof avoiding the use of 
metric regularity and Mordukhovich criterion.
\begin{proposition}Let $C,D$ be two closed subsets of $\R^n$ and 
$\bar x\in C\cap D$. Assume 
$$[-N_C(\bar{x})]\cap N_D(\bar{x})=\{0\}.$$ 
Then there exist a neighborhood $U\subset\R^n\times\R^n$ of $(\bar x,\bar x)$ 
and a positive constant $c$ for which
\begin{equation}\label{regular}
\dist(0,\partial L_{C,D}(x,y))\geq c\|x-y\|>0,
\end{equation}
whenever $(x,y)\in U\cap[0<L_{C,D}<+\infty]$. In other words $L_{C,D}$ has the 
Kurdyka-\L ojasiewicz property at 
$(\bar x,\bar x)$ with $\varphi(s)=\frac{1}{c}\sqrt{2s}$.
\end{proposition}
\noindent
{\bf Proof.} We argue by contradiction. Let $(x_k,y_k)\in [0<L<+\infty]$ be a 
sequence converging to $(\bar x,\bar x)$ such that 
$$\frac{\dist(0,\partial L_{C,D}(x_k,y_k))}{\|x_k-y_k\|}\leq\frac{1}{k+1}.$$
In view of (\ref{lojaset}), there exists 
$(u_k,v_k)\in N_{C}(x_k)\times N_D(y_k)$ such that 
$$\|x_k-y_k\|^{-1}[\|y_k-x_k-u_k\|+\|x_k-y_k-v_k\|]\leq\frac{\sqrt 2}{k+1}.$$
Let $d\in S^{n-1}$ be a cluster point of $\|x_k-y_k\|^{-1}(y_k-x_k).$ Since 
$\|x_k-y_k\|^{-1}[(y_k-x_k)-u_k]$ converges to zero, $d$ is also a limit-point 
of $\|x_k-y_k\|^{-1}u_k\in N_C(x_k)$. Due to the closedness property of $N_C$, 
we therefore obtain $d\in N_C(\bar x)$. Arguing similarly with $v_k$ we obtain 
 $d\in -N_{D}(\bar x)$. So $d$ satifies both $\|d\|=1$ and 
$-d\in[-N_C(\bar{x})]\cap N_D(\bar{x})$, a contradiction.$\hfill\Box$
\begin{remark}{\rm (a) Observe that (\ref{regular}) reads
$$\left(\dist^2(y-x,N_C(x))+\dist^2(x-y,N_D(y))\,\right)^{\frac{1}{2}}
   \geq c\|x-y\|.$$
(b) for related, but different, results see \cite{LLM,RL}.}
\end{remark}

\noindent
\subsection{Transverse Manifolds} As pointed out in \cite{mallew} a nice 
example of regular intersection is given by the smooth notion of 
transversality. Let $M$ be a smooth submanifold of $\R^n$. For each $x$ in 
$M$, $T_xM$ denotes the tangent space to $M$ at $x$. 
Two submanifolds of $\R^n$ are called transverse at $x\in M\cap N$ if they 
satisfy
$$T_xM+T_xN=\R^n.$$ 
Since the  normal cone to $M$ 
at $x$ is given by $N_M(x)=T_x M^{\perp},$
we have 
$$N_M(x)\cap N_N(x)=N_M(x)\cap[-N_N(x)]=\{0\},$$ whenever $M$ and $N$ are transverse at $x$. And therefore $L_{M,N}$ 
satisfies the Kurdyka-\L ojasiewicz inequality near $x$ with 
$\varphi(s)=c^{-1}\sqrt{s}$, $c>0$.

The constant $c$ can be estimated \cite[Theorem 18]{mallew}. 
\medskip

\noindent
\subsection{Tame feasibility}\label{S:feasibility}
Assume that $C,D$ are two globally analytic subsets of $\R^n$. The graph of 
$L_{C,D}$ is given by 
$$\graph L_{C,D}= G\cap (C\times D\times \R),$$
where $G$ is the graph of the polynomial function 
$(x,y)\rightarrow\frac{1}{2}||x-y||^2$. From the stability property of 
o-minimal structures (see Definition \ref{D:omin} (i)-(ii)), we deduce that 
$L_{C,D}$ is globally subanalytic.

Let $(\bar{x},\bar{y})\in C\times D$. Using the fact that the bifunction $L_{C,D}$ 
satisfies the \L ojasiewicz inequality with $\varphi(s)=s^{1-\theta}$ (where 
$\theta\in(0,1]$) we obtain 
\begin{equation}\label{dist_an}
\left(\dist^2(y-x,N_C(x))+\dist^2(x-y,N_D(y))\,\right)^\frac{1}{2}
\geq\left(\frac{1}{2}\|x-y\|^2\right)^\theta
\end{equation}
for all $(x,y)$ with $x\neq y$ in a neighborhood of $(\bar{x},\bar{y})$.

\begin{remark}{\rm (a) One of the most interesting features of this inequality 
is that it is satisfied for possibly tangent sets.\\
(b) Let us also observe that the above inequality is satisfied if $C,D$ are 
real-analytic submanifolds of $\R^n$. This is due to the fact that 
$B((\bar{x},\bar{y}),r)\cap [C\times D]$ is globally subanalytic for all 
$r>0$.
}
\end{remark}

If $C,D$ and the square norm function $||\cdot||^2$ are definable in the same 
 o-minimal structure, we of course have a similar result. Namely, 
\begin{equation}\label{dist_omin}
\varphi'\left(\frac{1}{2}\|x-y\|^2\right)
  \left(\dist^2(y-x,N_C(x))+\dist^2(x-y,N_D(y))\,\right)^{\frac{1}{2}}\geq 1,
\end{equation}
for all $(x,y)$ in a neighborhood of a critical point of $L$.

Although (\ref{dist_omin}) is simply a specialization of a result of 
\cite{BDLS}, this form of the Kurdyka-\L ojasiewicz inequality can be very 
useful in many contexts. 
 
\begin{remark}{\rm 
(a) Note that if $C,D$ are definable in the log-exp structure then (\ref{dist_omin}) holds. This is due to the fact that the square norm
  is also definable in this structure.\\
(b) If $C,D$ are semilinear sets, $L_{C,D}$ is not semilinear but it is however semialgebraic.\\
(c) Many concrete problems involving definable data and for which
 the projection operators are easily computable are given in \cite{mallew}.
}
\end{remark}


\begin{thebibliography}{99}

\bibitem{AMA2004}\textsc{Absil, P.-A., Mahony, R. , Andrews, B.,} Convergence 
of the iterates of descent methods for analytic cost functions, SIAM J. Optim.
{\bf 16} (2005), pp. 531--547.

\bibitem{AMS}\textsc{Absil, P.-A., Mahony, R. , Sepulchre, R.,} Optimization 
Algorithms on Matrix Manifolds, Princeton University Press, 2008.

\bibitem{ADG07}\textsc{Aragon, A., Dontchev, A. , Geoffroy, M.,} Convergence 
of the proximal point method for metrically regular mappings,   ESAIM Proc., {\bf 17}, EDP Sci., Les Ulis, (2007), pp 1--8.

\bibitem{attbol}{\sc Attouch, H., Bolte, J.,} On the convergence of the 
proximal algorithm for nonsmooth functions involving analytic features, Math. 
Program. Ser. B {\bf 116} (2009), pp. 5--16.

\bibitem{abrs}\textsc{Attouch, H., Bolte, J., Redont, P., Soubeyran, A.,}  
Alternating proximal algorithms for weakly coupled convex minimization 
problems.  Applications to dynamical games and PDE's, J. of Convex Analysis 
{\bf 15} (3) (2008), pp. 485--506.

\bibitem{ARS}{\sc  Attouch, H.,  Redont P.,  Soubeyran A.,} A new class of
alternating proximal minimization algorithms with costs to move, SIAM Journal
on Optimization {\bf 18} (2007), pp. 1061-1081.

\bibitem{AS}{\sc Attouch, H., Soubeyran A.,} Inertia and reactivity in 
decision making as cognitive variational inequalities, Journal of Convex 
Analysis {\bf 13} (2006), pp. 207--224.

%\bibitem{AS2}{\sc  Attouch, H., Soubeyran A.,} ``Worthwhile to move behaviours'' as temporary satisficing 
%processes with not too much sacrificing: how to overcome inertia and frictions?
% Journal of Mathematical Psychology, submitted, (2007). 

\bibitem{AE}{\sc Aubin, J.-P., Ekeland, I.,} {\it Applied Nonlinear Analysis}, 
Wiley-Interscience, 1984.

\bibitem{Aus}{\sc Auslender, A.,} {\it Optimisation, M\'ethodes Num\'eriques}, 
Masson, Paris, 1976.

\bibitem{BB}\textsc{Bauschke, H.H., Borwein, J.M.,} On projection algorithms 
for solving convex feasibility problems, SIAM review  {\bf 38} (3) (1996), 
 pp. 367-426.

\bibitem{BCN}{\sc H. H. Bauschke, P. L. Combettes and D. Noll}, {\em Joint 
minimization with alternating Bregman proximity operators}, Pacific Journal of 
Optimization, vol. 2, n$^\circ$ 3 (2006), pp.~401-424.

\bibitem{BR1990}\textsc{Benedetti, R., Risler, J.-J.,} \textit{Real Algebraic 
and Semialgebraic Sets}, Hermann, \'Editeur des Sciences et des Arts, Paris, 
 1990.

\bibitem{Ber} {\sc Bertsekas, D.,} {\it Nonlinear Optimisation}, Athena, Belmont, 
Massachusetts, second edition, 1999.

%\bibitem {BM88}\textsc{Bierstone, E., Milman, P.,} Semianalytic and
%subanalytic sets, \textit{IHES Publ. Math.} \textbf{67} (1988), 5--42.

\bibitem{BcR1998}\textsc{Bochnak, J., Coste, M., Roy, M.-F., }\textit{Real
Algebraic Geometry}, Springer, 1998.

\bibitem{BDLM}\textsc{Bolte, J., Daniilidis, A., Ley, O., Mazet, L.,} 
Characterizations of \L ojasiewicz inequalities and applications: subgradient flows, talweg, convexity,  Transactions 
of the AMS, 2009, to appear.

\bibitem{BDL04}\textsc{Bolte, J., Daniilidis, A. , Lewis, A.,} The 
\L ojasiewicz inequality for nonsmooth subanalytic functions with applications 
to subgradient dynamical systems, SIAM J. Optim. {\bf 17} (2006), 
pp. 1205--1223.

\bibitem {BDL04b}\textsc{Bolte, J., Daniilidis, A., Lewis, A.,} A nonsmooth 
Morse-Sard theorem for subanalytic functions, J. Math. Anal. Appl. {\bf 321} 
 (2006), pp. 729--740.

\bibitem {BDLS}\textsc{Bolte, J., Daniilidis, A., Lewis, A., Shiota, M.,} 
Clarke subgradients of stratifiable functions, SIAM J. Optim. {\bf 18} 
(2007), pp. 556--572.

\bibitem {BDE}\textsc{Bruckstein A., Donoho D., Elad M.,} 
 From sparse solutions of systems of equations to sparse modeling of signals 
and images,
SIAM Review  \textbf{51} (1) (2009), pp. 34--81. 

\bibitem {BC}\textsc{Byrne, C., Censor, Y.,} 
Proximity function minimization using multiple Bregman projections, with 
applications to split feasibility and Kullback-Leibler distance minimization,
Annals of Operations Research {\bf 105} (2001), pp. 77--98. 

\bibitem {CWB}\textsc{Cand\`es E., Wakin M., Boyd S.,} 
Enhancing sparsity by reweighting $\ell^1$ minimization, 
J. Fourier Anal. Appl. {\bf 14} (2008), pp. 877--905.

%\bibitem {CG}\textsc{Chen, Y., Gazzale, R.,} When does learning in games 
%generate convergence to Nash equilibria? The role of supermodularity in an
%experimental setting, American Economic Review, {\bf 19}, no. 5, (2004), 
%1505--1535.

%\bibitem {Chartrand}\textsc{Chartrand, R.,} 
%Exact reconstruction of sparse signals via nonconvex minimization,
%Signal Process. Lett. {\bf 14}(10), 707-710 (2007).

%\bibitem {SCY}\textsc{Saab R., Chartrand R., Yilmaz O.,} 
%Stable sparse approximation via nonconvex minimization, 
%in 33rd International Conference on Acoustics, Speech and Signal Processing 
%(ICASSP), 2008.
%
%\bibitem {clarke1983}\textsc{Clarke, F.H.,} \textit{Optimization and
% Nonsmooth Analysis}, Wiley-Interscience, (New York, 1983).

%\bibitem {cLSW98}\textsc{Clarke, F.H., Ledyaev, Yu., Stern, R.I. , Wolenski,
%P.R.,} \emph{Nonsmooth analysis and control theory}, Graduate texts in
%Mathematics \textbf{178}, Springer-Verlag, New-York, 1998.

\bibitem{Comb}\textsc{Combettes, P.L.,} Signal recovery by best feasible 
approximation, IEEE transactions on Image Processing {\bf 2} (1993), 
pp. 269--271.
%
%\bibitem{CoTr}\textsc{Combettes, P.L., Trussell, H.J.,} Method of successive 
%projections for finding a common point of sets in metric spaces, Journal of 
%Optimization Theory and Applications, {\bf 67} (3), (1990), 487-507.

\bibitem{CP04} {\sc Combettes, P. L., Pennanen, T.,} Proximal methods for cohypomonotone operators.  SIAM J. Control Optim.  {\bf 43}  (2004),  pp. 731--742


\bibitem{CombWajs}\textsc{Combettes, P.L., Wajs, V.} 
Signal recovery by proximal forward-backward splitting,  
Multiscale Model. Simul. {\bf 4}  (2005), pp. 1168--1200.
%
%\bibitem {coste99}\textsc{Coste, M.,} \textit{An introduction to o-minimal
%geometry}, RAAG Notes, 81 p., Institut de Recherche Math\'{e}matiques de
%Rennes, November 1999.

\bibitem{CT} {\sc Csiszar, I., Tusnady, G.,} Information geometry and 
alternating minimization procedures, Statistics and Decisions (supplement 1)
(1984), pp. 205--237.

\bibitem{Denk} {\sc Denkowska, Z., Stasica, J.,} Ensembles sous-analytiques 
\`a la Polonaise, Travaux en cours {\bf 69}, Hermann, Paris 2008, to appear.
%
%\bibitem{De}\textsc{Deutsch, F.,} Best approximation in inner product spaces, 
%Springer, New York, 2001.

\bibitem{donoho}{\sc Donoho, D. L.}, 
Compressed Sensing, 
IEEE Trans. Inform. Theory {\bf 4} (2006), pp. 1289--1306. 

\bibitem {Dries}\textsc{van den Dries, L.,} Tame topology and o-minimal 
structures, London Mathematical Society Lecture Note Series, {\bf 248}, 
Cambridge University Press, Cambridge, (1998) x+180 pp.

\bibitem{Edel} {\sc Edelman, A., Arias,  A., Smith, S.T.,} The geometry of algorithms with orthogonality constraints,  SIAM J. Matrix Anal. Appl.  {\bf 20}   (2) (1999),  pp. 303--353.

\bibitem{Hare} {\sc Hare, W.,  Sagastiz\'abal, C.,} Computing proximal points of nonconvex functions,  Math. Program. Ser. B   {\bf 116}  (2009), pp. 221--258.

%\bibitem{Ha}\textsc{Haraux, A.,} A hyperbolic variant of Simon's convergence 
%theorem, Lectures Notes in Pure and Appl. Math. \textbf{215}, Dekker-New York 
%(2001), 255-264.

\bibitem{gabrielov} {\sc Gabrielov, A.,} Complements of subanalytic sets and 
existential formulas for analytic functions, Inventiones Math. \textbf{125} 
(1996), pp. 1--12.
 
\bibitem{GrubP} {\sc Grubisi\'c I., Pietersz R.},  Efficient rank reduction of correlation matrices, Linear Algebra and its Applications {\bf 422} (2007), pp. 629--653.

%\bibitem{Higham} {\sc Highmam, N.,} Computing a nearest symmetric correlation 
%matrix; a problem from finance. IMA Journal of Numerical Analysis, {\bf 22}, 
%no. 3, (2002), 329-343.

\bibitem{Ioffe} {\sc Ioffe, A.D.,} Metric regularity and subdifferential 
calculus, (Russian)  Uspekhi Mat. Nauk  \textbf{55} (2000), no. 3 (333), 
pp.  103--162; translation in Russian Math. Surveys \textbf{55} (2000), pp.  
501--558.

\bibitem{IPS} {\sc Iusem, A. N., Pennanen, T., Svaiter, B. F.} Inexact variants of the proximal point algorithm without monotonicity.  SIAM J. Optim.  {\bf 13}  (2003),  no. 4, pp. 1080--1097.


\bibitem{Kurdyka98}\textsc{Kurdyka, K.,} On gradients of functions definable 
in o-minimal structures, Ann. Inst. Fourier  \textbf{48} (1998), pp. 769--783.

\bibitem{mallew}{\sc Lewis, A.S., Malick, J.,} Alternating projection on 
manifolds,  Math. Oper. Res.  {\bf 33}  (2008),  pp. 216--234.

\bibitem{LLM}{\sc Lewis, A.S., Luke, D.R., Malick, J.,} Local Linear 
Convergence for Alternating and Averaged Nonconvex Projections, Found. Comput. 
Math. {\bf 9} no. 4 (2009), pp. 485-513.

\bibitem{PL.Lions} {\sc Lions, P.L.,} On the Schwarz alternating method. III. 
A variant for nonoverlapping subdomains  in Third International Symposium on 
Domain Decomposition Methods for Partial Differential Equations, T. F. Chan, 
R. Glowinski, J. P\'eriaux and O. Widlund, eds., SIAM (1990), pp. 202--231.

\bibitem{Loja63}\textsc{\L ojasiewicz, S.,} Une propri\'et\'e topologique des 
sous-ensembles analytiques r\'eels, in: \textit{Les \'Equations aux 
D\'eriv\'ees Partielles}, pp. 87--89, \'Editions du centre National de la 
Recherche Scientifique, Paris 1963.

\bibitem{Loja93}\textsc{\L ojasiewicz, S.,} Sur la g\'eom\'etrie semi- et
sous-analytique, Ann. Inst. Fourier \textbf{43} (1993), pp. 1575--1595.

%\bibitem {Marker2002}\textsc{Marker, D., }\textit{Model theory. An
%introduction,} Graduate Texts in Mathematics \textbf{217}, (Springer, %2002).

%\bibitem {Milman02}\textsc{Milman, E.,} The semi-algebraic theory of stochastic
%games, Math.\, Oper.\, Res. \textbf {27} (2002), 401--418.

%\bibitem {MS}\textsc{Monderer, D., Shapley, L.,} Fictitious play property  for 
%games with identical players, Journal of Economic Theory, {\bf 68}, (1996), 
%258-265.

%\bibitem{Morduk76} {\sc  Mordukhovich, B.,}  Maximum principle in the problem 
%of time optimal response with nonsmooth constraints,  J. Appl. Math. 
%Mech., {\bf 40} (1976), 960--969 ; [translated from Prikl. Mat. Meh. 40 
%(1976), 1014--1023] 

\bibitem{RL}{\sc Luke R. D.,} Finding best approximation pairs relative 
to a convex and prox-regular set in a Hilbert space, SIAM J. Optim. Vol. 19, 
No. 2, pp. 714-739.

\bibitem{Morduk}{\sc Mordukhovich, B.}, {\it Variational analysis and 
generalized differentiation. I. Basic theory}, Grund\-lehren der 
Mathematischen Wissenschaften, {\bf 330}, Springer-Verlag, Berlin, 2006, 
xxii+579 pp.
%
%\bibitem{NestPol}{\sc Nesterov, Yu.},  Accelerating the cubic regularization 
%of Newton's method on convex problems, Math. Program. {\bf 112} (2008), no. 1, 
%Ser. B, 159--181.

\bibitem{Nesterov}{\sc Nesterov, Yu.}, 
{\it Introductory lectures on convex optimization. A basic course.} 
Applied Optimization, {\bf 87}, Kluwer Academic Publishers, Boston, MA, 2004. 
xviii+236 pp. 

\bibitem{Johannis}{\sc von Neumann, J.}, 
{\em Functional Operators}, 
Annals of Mathematics studies {\bf 22}, Princeton University Press, 1950. 

\bibitem {Rock98}\textsc{Rockafellar, R.T. , Wets, R.,} \textit{Variational
Analysis}, Grundlehren der Mathematischen Wissenschaften,  \textbf{317},
Springer, 1998.


%\bibitem{Sim83}\textsc{Simon, L., } Asymptotics for a class of non-linear
%evolution equations, with applications to geometric problems, Ann.  Math. 
%\textbf{118} (1983), 525-571.

\bibitem {Shiota}\textsc{Shiota, M.,} \textit{Geometry of subanalytic and
semialgebraic sets}, Progress in Mathematics \textbf{150}, Birkh\"{a}user, 
Boston, 1997.

\bibitem{Tseng} {\sc Tseng, P.}, 
Convergence of a block coordinate descent method for nondifferentiable minimization,  
J. Optim. Theory Appl.  {\bf 109}  (2001),  pp. 475--494.

%\bibitem{Smale}\textsc{Smale, S., } Mathematical problems for the next 
%century, Mathematical Intelligence, {\bf 20}, no. 2, (1998), 7-15.

\bibitem {WW}\textsc{Widrow, B., Wallach, E.,} \textit{Adaptive Inverse 
Control}, Prentice-Hall, New Jersey, 1996.

\bibitem{Wilkie}{\sc Wilkie, A. J.,}  Model completeness results for 
expansions of the ordered field of real numbers by restricted Pfaffian 
functions and the exponential function, J. Amer. Math. Soc. {\bf 9}  
(1996), pp. 1051--1094.

\end{thebibliography}
\end{document}